\documentclass[a4paper,oneside,10pt]{amsart}
\usepackage{amssymb}
\newtheorem{thm}{Theorem}[section]
\newtheorem{lem}[thm]{Lemma}
\newtheorem{cor}[thm]{Corollary}
\newtheorem{prop}[thm]{Proposition}
\newtheorem{rem}[thm]{Remark}
\newtheorem{exmp}[thm]{Example}
\newtheorem{defn}[thm]{Definition}
\title{On the direct indecomposability of infinite irreducible Coxeter groups and the Isomorphism Problem of Coxeter groups}
\author{Koji NUIDA}
\thanks{The author is supported by JSPS Research Fellowship (No.\ 16-10825)}
\begin{document}
\maketitle
\begin{abstract}
In this paper we prove that any irreducible Coxeter group of infinite order is directly indecomposable as an abstract group, without the finite rank assumption.
The key ingredient of the proof is that we can determine, for an irreducible Coxeter group $W$, the centralizers in $W$ of the normal subgroups of $W$ that are generated by involutions.
As a consequence, we show that the problem of deciding whether two general Coxeter groups are isomorphic, as abstract groups, is reduced to the case of irreducible Coxeter groups, without assuming the finiteness of the number of the irreducible components or their ranks.
We also give a description of the automorphism group of a general Coxeter group in terms of those of its irreducible components.
\end{abstract}
\section{Introduction}
In this paper, we prove that all infinite irreducible Coxeter groups are directly indecomposable as abstract groups (Theorem \ref{thm:indecomposabilityofW}).\\ \indent
Regarding direct indecomposability of Coxeter groups, it is well known that there exist finite irreducible Coxeter groups which are directly decomposable (such as the Weyl group $G_2$).
On the other hand, for infinite irreducible Coxeter groups, no general result has been known until recently.
In a recent paper \cite{Par}, L.~Paris proved the direct indecomposability of all infinite irreducible Coxeter groups of finite rank, by using certain special elements called essential elements which are examined also in \cite{Kra}.
However, by definition, a Coxeter group of infinite rank never possesses an essential element, so that the proof cannot be applied directly to the case of infinite rank.\\ \indent
Our result here is obtained by a different approach.
Let $W$ be an irreducible Coxeter group whose order is infinite, possibly of infinite rank.
We give a complete description of the centralizer $C$ of any normal subgroup $N$ of $W$ which are generated by involutions (Theorem \ref{thm:Zofinvolutivenormalsubgroup}).
From the description it follows that, unless $N=\{1\}$ or $C=\{1\}$, there is a subgroup $H\subsetneq W$ which contains both $N$ and $C$.
Once this is proved, the direct indecomposability of $W$ is clear, since any direct factor of $W$ is a normal subgroup and is generated by involutions (since it is a quotient of $W$), and its centralizer contains the complementary factor.\\ \indent
As a consequence of the direct indecomposability of infinite irreducible Coxeter groups, we give results on the isomorphisms between two Coxeter groups (Theorem \ref{thm:isomorphismproblemforW}).
Since we also know how each finite irreducible Coxeter group decomposes into directly indecomposable factors, our results imply that we can determine whether or not two given Coxeter groups are isomorphic if we can determine which infinite irreducible Coxeter groups are isomorphic.
In addition, our results also give certain decompositions of an automorphism of a general Coxeter group $W$ (Theorem \ref{thm:automorphismgroup}).
One decomposition describes its form from the viewpoint of the directly indecomposable decomposition of $W$; another decomposition describes its form from the viewpoint of the decomposition $W=W_{\rm fin}\times W_{\rm inf}$, where $W_{\rm fin}$ (resp.\ $W_{\rm inf}$) is the product of the finite (resp.\ infinite) irreducible components of $W$ in the given Coxeter system.
Note that these results can also be deduced from the Krull-Remak-Schmidt Theorem in group theory, if the Coxeter group has a composition series.
Theorem \ref{thm:isomorphismproblemforW} is also a generalization of Theorem 2.1 of \cite{Par}; our proof here is similar to, but slightly more delicate than that in \cite{Par}, by the lack of finiteness of the ranks.
Note also that, in another recent paper \cite{Mih-Rat-Tsc}, M.~Mihalik, J.~Ratcliffe and S.~Tschantz also examined the ``Isomorphism Problem'' (namely, the problem of deciding which Coxeter groups are isomorphic) for the case of finite ranks, by a highly different approach.\\ \indent
{\bf Contents.}\ 
Section \ref{sec:preliminaries} collects the preliminary facts and results.
In Section \ref{sec:generalgroups}, we give some remarks on general groups, especially on the definition and properties of the core subgroups.
Sections \ref{sec:Coxetergraph} and \ref{sec:rootsystem} summarize definitions, notations and properties of Coxeter systems, Coxeter graphs and root systems of Coxeter groups.
In Section \ref{sec:decompositionofw_0}, we recall a method, given by V.~Deodhar \cite{Deo}, for decomposing the longest element of any finite parabolic subgroup into pairwise commuting reflections.
Owing to this decomposition, we can compute easily the action of the longest element on a root, even if it is not contained in the root system of the parabolic subgroup.
As an application, in Section \ref{sec:decompositionoffiniteW}, we determines all irreducible Coxeter groups of which the center is a nontrivial direct factor.
(This is not a new result, but is included there since the result is used in the following sections.)
Some properties of normalizers of parabolic subgroups are summarized as Section \ref{sec:noteonnormalizer}.\\ \indent
Our main results are stated and proved in Section \ref{sec:mainresults}.
The direct indecomposability of infinite irreducible Coxeter groups is shown in Section \ref{sec:indecomposability} (Theorem \ref{thm:indecomposabilityofW}).
Note that the theorem also determines all nontrivial direct product decompositions of finite irreducible Coxeter groups.
In Section \ref{sec:isomorphismproblem}, we reduce the Isomorphism Problem of general Coxeter groups to the case of infinite irreducible ones (Theorem \ref{thm:isomorphismproblemforW}).
In the proof, we consider such a problem in a slightly wider context (Theorem \ref{thm:isomorphismproblem}) and then our result is deduced.
Moreover, another result in Section \ref{sec:automorphismgroup} describes the automorphism group of a general Coxeter group in terms of those of the irreducible components (Theorem \ref{thm:automorphismgroup} (ii)).
Note that a Coxeter group possesses some `natural' automorphisms, which map each irreducible component onto a component isomorphic to the original one.
We also give a characterization of Coxeter groups for which the group of the `natural' automorphisms has finite index in the whole automorphism group (Theorem \ref{thm:automorphismgroup} (iii)).\\ \indent
Our proof of Theorem \ref{thm:indecomposabilityofW} is based on our description of the centralizers of the normal subgroups, which are generated by involutions, in irreducible Coxeter groups (Theorem \ref{thm:Zofinvolutivenormalsubgroup}).
This theorem is proved in Section \ref{sec:proofoftheorem_Sec3}, by using a description (given in Sections \ref{sec:lemmaforcores}--\ref{sec:coreofinfinitecase}) of core subgroups of normalizers of parabolic subgroups.\\ \indent
{\bf Acknowledgement.}\quad 
I would like to express my deep gratitude especially to Itaru Terada and Kazuhiko Koike, for their precious advice and encouragement.
\section{Preliminaries}\label{sec:preliminaries}
\subsection{Notes on general groups}\label{sec:generalgroups}
In this paper, we treat two kinds of direct products of groups $G_\lambda$ with (possibly infinite) index set $\Lambda$; the {\it complete} direct product (whose elements $(g_\lambda)_\lambda$ are all the maps $\Lambda \to \bigsqcup_{\mu \in \Lambda}G_\mu$, $\lambda \mapsto g_\lambda$ such that $g_\lambda \in G_\lambda$) and the {\it restricted} direct product (consisting of all the elements $(g_\lambda)_\lambda$ such that $g_\lambda$ is the unit element of $G_\lambda$ for all but finitely many $\lambda \in \Lambda$).
Note that these two products coincide if $|\Lambda|<\infty$.
Since here we treat mainly the latter type rather than the former one, we let {\it the term ``direct product'' alone and the symbol $\prod$ mean the restricted direct product} throughout this paper.
(The complete one also appears in this paper, always together with notification.)\\ \indent
For two groups $G$, $G'$, let ${\rm Hom}(G,G')$, ${\rm Isom}(G,G')$ denote the sets of all homomorphisms, isomorphisms $G \to G'$ respectively.
Put ${\rm End}(G)={\rm Hom}(G,G)$ and ${\rm Aut}(G)={\rm Isom}(G,G)$.
The following lemma is easy, but will be referred later.
\begin{lem}
\label{lem:conditionforvanishingcenter}
Assume that the center $Z(G)$ of a group $G$ is either trivial or a cyclic group of prime order.
Then the following three conditions are equivalent:\\
{\bf (I)} $Z(G)=1$ or $Z(G)$ is not a direct factor of $G$.\\
{\bf (II)} If $f \in {\rm Hom}(G,Z(G))$, then $f(Z(G))=1$.\\
{\bf (III)} If $G'$ is a direct product of (arbitrarily many) cyclic groups of prime order and $f \in {\rm Hom}(G,G')$, then $f(Z(G))=1$.
\end{lem}
\begin{proof}
This is trivial if $Z(G)=1$, so that we assume that $Z(G)$ is a cyclic group of prime order.
Note that the implication (III) $\Rightarrow$ (II) is obvious.\\
{\bf (I) $\Leftrightarrow$ (II):}\ 
If (I) is not satisfied, and $G=Z(G) \times H$, then the projection $G \to Z(G)$ does not satisfy the conclusion of (II).
Conversely, if $f \in {\rm Hom}(G,Z(G))$ and $f(Z(G)) \neq 1$, then $f(Z(G))=Z(G)$, $\ker f \cap Z(G)=1$ (since $Z(G)$ is simple) and so we have $G=Z(G) \times \ker f$.\\
{\bf (II) $\Rightarrow$ (III):}\ 
This is clear if $G'$ itself is a cyclic group of prime order (by noting that ${\rm Hom}(\mathbb{Z}/p\mathbb{Z},\mathbb{Z}/\ell\mathbb{Z})=1$ for distinct primes $p$, $\ell$).
For a general case, apply it to the composite map $\pi \circ f$ for every projection $\pi$ from $G'$ to one of its factors.
\end{proof}
Here we define the following multiplication for the set ${\rm Hom}(G,Z(G))$ by which it forms a monoid.
First, we define a map ${\rm Hom}(G,Z(G)) \to {\rm End}(G)$, $f \mapsto f^\flat$ by
\begin{equation*}
f^\flat(w)=wf(w)^{-1} \text{ for all } w \in G.
\end{equation*}
This is well defined since $Z(G)$ is abelian.
The image of $H \subset {\rm Hom}(G,Z(G))$ by the map is denoted by $H^\flat$.
Now define the product $f*g$ of two elements $f,g \in {\rm Hom}(G,Z(G))$ by
\begin{equation*}
(f*g)(w)=f(w)g(w)f \circ g(w)^{-1} \text{ for all } w \in G.
\end{equation*}
This is also well defined, and then ${\rm Hom}(G,Z(G))$ forms a monoid with the trivial map (denoted by $1$) as the unit element (for example, we have the associativity
\begin{equation}
\label{eq:associativity}
\begin{split}
\bigl((f*g)*h\bigr)(w)&=\bigl(f*(g*h)\bigr)(w)\\
&=f(w)g(w)h(w)f \circ g(w)^{-1}f \circ h(w)^{-1}g \circ h(w)^{-1}f \circ g \circ h(w)\\
&=(f*h)(w)f^\flat \circ g \circ h^\flat(w)
\end{split}
\end{equation}
for $f,g,h \in {\rm Hom}(G,Z(G))$).
Let ${\rm Hom}(G,Z(G))^\times$ denote the group of invertible elements of ${\rm Hom}(G,Z(G))$ with respect to the multiplication $*$.
On the other hand, ${\rm End}(G)$ also forms a monoid with composition of maps as multiplication; then the group of invertible elements in the monoid ${\rm End}(G)$ is precisely the group ${\rm Aut}(G)$.\\ \indent
Moreover, the group ${\rm Aut}(G)$ acts on the monoids ${\rm Hom}(G,Z(G))$ and ${\rm End}(G)$ by
\begin{equation*}
h \cdot f=h \circ f \circ h^{-1}\ \text{ for } h \in {\rm Aut}(G),\ f \in {\rm Hom}(G,Z(G)) \text{ or } {\rm End}(G).
\end{equation*}
\begin{lem}
\label{lem:HomtoEndisembedding}
{\bf (i)} The map $f \mapsto f^\flat$ is an injective homomorphism ${\rm Hom}(G,Z(G)) \to {\rm End}(G)$ of monoids compatible with the action of ${\rm Aut}(G)$.\\
{\bf (ii)} For $f \in {\rm Hom}(G,Z(G))$, the following three conditions are equivalent:\\ \indent
{\bf (I)} $f \in {\rm Hom}(G,Z(G))^\times$.\qquad 
{\bf (II)} $f^\flat \in {\rm Aut}(G)$.\\ \indent
{\bf (III)} The restriction $f^\flat|_{Z(G)}$ is an automorphism of $Z(G)$.\\
{\bf (iii)} If $H \subset {\rm Hom}(G,Z(G))^\times$ is a subgroup invariant under the action of ${\rm Aut}(G)$, then its image $H^\flat$ is a normal subgroup of ${\rm Aut}(G)$.
\end{lem}
\begin{proof}
The claim (i) is straightforward, while (iii) follows from (i), (ii) and definition of the action of ${\rm Aut}(G)$.
From now, we prove (ii).
The implication (I) $\Rightarrow$ (II) is obvious.
On the other hand, (II) implies (III) since any automorphism preserves the center.
Moreover, if (III) is satisfied, then we can construct the inverse element $f'$ of $f \in {\rm Hom}(G,Z(G))$ by $f'(w)=(f^\flat|_{Z(G)})^{-1}\bigl(f(w)\bigr)^{-1}$ ($w \in G$); we have
\begin{equation*}
\begin{split}
(f'*f)(w)=f'(w)f(w)f'(f(w))^{-1}&=f'(wf(w)^{-1})f(w)\\
&=(f^\flat|_{Z(G)})^{-1}\bigl(f(wf(w)^{-1})\bigr)^{-1}f(w)\\
&=(f^\flat|_{Z(G)})^{-1}\bigl(f^\flat(f(w))\bigr)^{-1}f(w)\\
&=f(w)^{-1}f(w)=1,
\end{split}
\end{equation*}
so that $f'*f=1$.
Similarly, we have $f*f'=1$.
Hence the claim holds.
\end{proof}
\begin{lem}
\label{lem:Homofabelian}
If a group $G$ is abelian, then the embedding ${\rm Hom}(G,Z(G)) \to {\rm End}(G)$, $f \mapsto f^\flat$, is an isomorphism with inverse map $f \mapsto f^\flat$.
Moreover, its restriction is an isomorphism ${\rm Hom}(G,Z(G))^\times \to {\rm Aut}(G)$.
\end{lem}
\begin{proof}
Note that $Z(G)=G$, so that ${\rm Hom}(G,Z(G))={\rm End}(G)$ as sets.
Thus the map ${\rm End}(G) \to {\rm Hom}(G,Z(G))$, $f \mapsto f^\flat$ is well defined.
Now we have $(f^\flat)^\flat(w)=wf^\flat(w)^{-1}=f(w)$ for all $f \in {\rm End}(G)$ and $w \in G$, so that $(f^\flat)^\flat=f$.
Thus the first claim holds.
Now the second one follows from Lemma \ref{lem:HomtoEndisembedding} (ii).
\end{proof}
Note that, if $G=G_1 \times G_2$, then the sets ${\rm Hom}(G_i,Z(G))$ ($i=1,2$) are embedded into ${\rm Hom}(G,Z(G))$ via the map $f \mapsto f \circ \pi_i$ (where $\pi_i$ is the projection $G \to G_i$).
Each ${\rm Hom}(G_i,Z(G))$ forms a submonoid of ${\rm Hom}(G,Z(G))$.
Moreover, the above formula of the inverse element $f'$ of $f \in {\rm Hom}(G,Z(G))$ implies that, $f \in {\rm Hom}(G_i,Z(G))$ is invertible in ${\rm Hom}(G_i,Z(G))$ if and only if it is invertible in ${\rm Hom}(G,Z(G))$.
Thus the notation ${\rm Hom}(G_i,Z(G))^\times$ is unambiguous.
\begin{lem}
\label{lem:subgroupofHom}
{\bf (i)} Let $f,g \in {\rm Hom}(G,Z(G))$ such that $f(Z(G))=g(Z(G))=1$.
Then $f,g \in {\rm Hom}(G,Z(G))^\times$ and $(f*g)(w)=f(w)g(w)$ for all $w \in G$ (so that $f*g=g*f$ by symmetry).
Moreover, the map $w \mapsto f(w)^{-1}$ is the inverse element of $f$ in ${\rm Hom}(G,Z(G))^\times$.\\
{\bf (ii)} Suppose that $G=G_1 \times G_2$ and $Z(G_2)=1$.
Then ${\rm Hom}(G,Z(G))^\times=H_1 \rtimes H_2$ where $H_1={\rm Hom}(G_2,Z(G))$, $H_2={\rm Hom}(G_1,Z(G_1))^\times$.
Moreover, $H_1$ is abelian, $(f*g)(w)=f(w)g(w)$ for $f,g \in H_1$ and $f*g*f'=f^\flat \circ g \circ (f^\flat)^{-1}$ for $f \in H_2$ and $g \in H_1$, where $f'$ is the inverse element of $f \in H_2$.
\end{lem}
\begin{proof}
{\bf (i)} By the hypothesis, $f^\flat$ is identity on $Z(G)$, so that $f$ is invertible by Lemma \ref{lem:HomtoEndisembedding} (ii) (and $g$ is so).
The other claims follow from definition (note that now $f \circ g=1$).\\
{\bf (ii)} Note that $Z(G)=Z(G_1)$ by the hypothesis.
Then by (i), $H_1$ is an abelian subgroup of ${\rm Hom}(G,Z(G))^\times$ in which the multiplication is as in the statement.\\ \indent
For $f \in H_2$ and $g \in H_1$, the formula (\ref{eq:associativity}) implies that $f*g*f'$ is as in the statement (note that $f*f'=1$ and ${f'}^\flat=(f^\flat)^{-1}$).
In particular, we have $f*g*f'(G_1) \subset f^\flat \circ g(G_1)=1$, since $f' \in H_2$ and so ${f'}^\flat(G_1) \subset G_1$.
This means that $f*g*f' \in H_1$.
Since obviously $H_1 \cap H_2=1$, we have $H_1H_2=H_1 \rtimes H_2$.\\ \indent
Finally, let $f \in {\rm Hom}(G,Z(G))^\times$.
Take $g \in H_1$ such that $g(w)=f \circ \pi_2(w)^{-1}$ where $\pi_2$ is the projection $G \to G_2$ (this is the inverse element of $f \circ \pi_2 \in H_1$).
Then for $w \in G_2$, we have
\begin{equation*}
(g*f)(w)=g(w)f(w)g\bigl(f(w)\bigr)^{-1}=f(w)^{-1}f(w)=1
\end{equation*} 
since $g(Z(G))=1$.
This means that $g*f \in {\rm Hom}(G_1,Z(G_1))$, while it is invertible since both $f$ and $g$ are so.
Thus we have $g*f \in H_2$ and $f=(f \circ \pi_2)*g*f \in H_1H_2$.
Hence ${\rm Hom}(G,Z(G))^\times=H_1 \rtimes H_2$.
\end{proof}
In the proof of our results, we use the following notion.
For a group $G$, we write $H \leq G$, $H \lhd G$ if $H$ is a subgroup, normal subgroup of $G$, respectively.
\begin{defn}
\label{defn:coresubgroup}
For $H \leq G$, define the {\it core} ${\rm Core}_G(H)$ of $H$ in $G$ to be the unique maximal normal subgroup of $G$ contained in $H$ (namely, $\bigcap_{w \in G}wHw^{-1}$).
\end{defn}
The following properties are deduced immediately from definition:
\begin{align}
\label{eq:coreisincreasing}
&\text{If } H_1 \leq H_2 \leq G, \text{ then } {\rm Core}_G(H_1) \subset {\rm Core}_G(H_2).\\
\label{eq:coreinsubgroup}
&\text{If } {\rm Core}_G(H) \leq H_1 \leq G, \text{ then } {\rm Core}_G(H) \subset {\rm Core}_G(H_1).\\
\label{eq:coreofintersection}
&\text{If } H_\lambda \leq G \text{ ($\lambda \in \Lambda$), then } {\rm Core}_G(\bigcap_{\lambda \in \Lambda}H_\lambda)=\bigcap_{\lambda \in \Lambda}{\rm Core}_G(H_\lambda).\\
\label{eq:subgroupoutsidecore}
&\text{If } H_1 \leq H_2 \leq G,\ w \in G \text{ and } wH_1w^{-1} \cap H_2=1, \text{ then } H_1 \cap {\rm Core}_G(H_2)=1.
\end{align}
\begin{lem}
\label{lem:coreoflimit}
Let $G_1 \leq G_2 \leq \dotsm$, $H_1 \leq H_2 \leq \dotsm$ be two infinite chains of subgroups of the same group such that $G_i \cap H_j=H_i$ for all $i \leq j$.
Put $G=\bigcup_{i=1}^{\infty}G_i$ and $H=\bigcup_{i=1}^{\infty}H_i$.
Then ${\rm Core}_G(H) \subset \bigcup_{i=1}^{\infty}{\rm Core}_{G_i}(H_i)$.
\end{lem}
\begin{proof}
It is enough to show that ${\rm Core}_G(H) \cap H_i \subset {\rm Core}_{G_i}(H_i)$ (or more strongly, ${\rm Core}_G(H) \cap H_i \lhd G_i$) for all $i$.
Note that the hypothesis implies $G_i \cap H=H_i$.
Then for $g \in G_i$ and $h \in {\rm Core}_G(H) \cap H_i$, we have $ghg^{-1} \in {\rm Core}_G(H)$ and $ghg^{-1} \in G_i$, so that $ghg^{-1} \in G_i \cap H=H_i$.
Thus the claim holds.
\end{proof}
The next lemma describes the centralizers of normal subgroups in terms of the cores of certain subgroups.
Before stating this, note the following easy facts:
\begin{align}
\label{eq:Zisnormal}
&\text{If } H \lhd G, \text{ then the centralizer } Z_G(H) \text{ of } H \text{ is also normal in } G.\\
\label{eq:Zreverseinclusion}
&\text{If } X_1,X_2 \subset G \text{ are subsets and } X_1 \subset Z_G(X_2), \text{ then } X_2 \subset Z_G(X_1).
\end{align}
\begin{lem}
\label{lem:Zandcore}
Let $H$ be the smallest normal subgroup of $G$ containing a subset $X \subset G$.
Then $Z_G(H)={\rm Core}_G(Z_G(X))=\bigcap_{x \in X}{\rm Core}_G(Z_G(x))$.
\end{lem}
\begin{proof}
The second equality follows from (\ref{eq:coreofintersection}).
For the first one, the inclusion $\subset$ is deduced from (\ref{eq:Zisnormal}) (since $Z_G(H) \subset Z_G(X)$).
For the other inclusion, the centralizer of ${\rm Core}_G(Z_G(X))$ in $G$ is normal in $G$ (by (\ref{eq:Zisnormal})) and contains $X$, so that it also contains $H$.
Thus the claim follows from (\ref{eq:Zreverseinclusion}).
\end{proof}
%
\subsection{Coxeter groups and Coxeter graphs}\label{sec:Coxetergraph}
Here we refer to \cite{Hum} for basic definitions and properties.
A pair $(W,S)$ of a group $W$ and its generating set $S$ is a {\it Coxeter system} (and $W$ itself is a {\it Coxeter group}) if $W$ has the presentation
\begin{equation*}
W=\langle S \mid (st)^{m(s,t)}=1 \text{ if } s,t \in S \text{ and } m(s,t)<\infty \rangle
\end{equation*}
where $m:S \times S \to \{1,2,\dotsc\} \cup \{\infty\}$ is a symmetric map such that $m(s,t)=1$ if and only if $s=t$.
$(W,S)$ is said to be {\it finite} ({\it infinite}) if the group $W$ is finite (infinite, respectively).
The cardinality of $S$ is called the {\it rank} of $(W,S)$ (or even of $W$).
Throughout this paper, {\it we do not assume, unless otherwise noticed, that the rank of $(W,S)$ is finite (or even countable)}.
Note that, owing to the well-known fact that the element $st \in W$ above has precisely order $m(s,t)$ in $W$, this map $m$ can be recovered uniquely from the Coxeter system $(W,S)$.\\ \indent
Two Coxeter systems $(W,S)$ and $(W',S')$ are said to be {\it isomorphic} if there is some $f \in {\rm Isom}(W,W')$ such that $f(S)=S'$.
Then there is a one-to-one correspondence (up to isomorphism) between Coxeter systems and the {\it Coxeter graphs}; which are simple(, loopless), undirected, edge-labelled graphs with labels in $\{3,4,\dotsc\} \cup \{\infty\}$.
The Coxeter graph $\varGamma$ corresponding to $(W,S)$ has the vertex set $S$, and two vertices $s,t \in S$ are joined in $\varGamma$ by an edge with label $m(s,t)$ if and only if $m(s,t) \geq 3$ (by convention, the labels `$3$' are usually omitted).
$\varGamma$ (or $(W,S)$) is said to be of {\it finite type} if $W$ is finite.
It is also well known that a full subgraph $\varGamma_I$ of $\varGamma$ with vertex set $I \subset S$ corresponds to a parabolic subgroup $W_I$ of $W$ generated by $I$ (or more precisely, to a Coxeter system $(W_I,I)$).\\ \indent
A Coxeter system $(W,S)$ is called {\it irreducible} if the corresponding Coxeter graph $\varGamma$ is connected.
In this case, $W$ is also said to be {\it irreducible}.
As is well known, $W$ is decomposed as the direct product of its {\it irreducible components}, which are the parabolic subgroups $W_I$ of $W$ corresponding to the connected components $\varGamma_I$ of $\varGamma$ (in this case, each subset $I$ is also said to be an {\it irreducible component} of $S$).
A parabolic subgroup $W_I \subset W$ is said to be {\it irreducible} if the Coxeter system $(W_I,I)$ is irreducible.
As we mentioned in Introduction, an irreducible Coxeter group may be directly decomposable (as an abstract group) in general.
Our main result determines which irreducible Coxeter group is indeed directly indecomposable.\\ \indent
In this paper, we use the following notations for some Coxeter graphs.
\begin{defn}
\label{defn:Coxetergraphs}
We use the notations in Fig.~1.
For each of the Coxeter graphs, let $s_i$ denote the vertex having label $i$.
Moreover, for each Coxeter graph $\varGamma(\mathcal{T}_n)$ in Fig.~1 ($\mathcal{T}=A$, $B$, $D$, $E$, $F$, $H$), let $\varGamma(\mathcal{T}_k)$ ($k<n$) be the full subgraph of $\varGamma(\mathcal{T}_n)$ on vertex set $\{s_i \mid 1 \leq i \leq k\}$.
For any $\mathcal{T}$, let $(W(\mathcal{T}),S(\mathcal{T}))$ be the Coxeter system corresponding to the Coxeter graph $\varGamma(\mathcal{T})$.
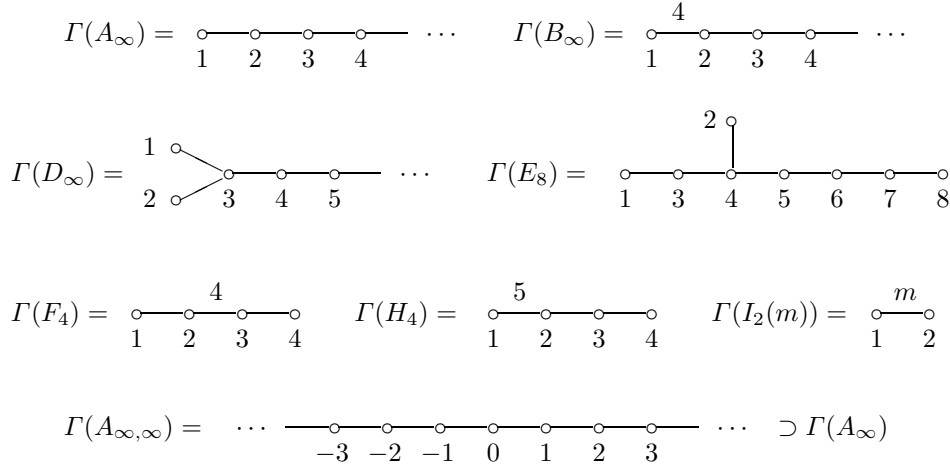
\begin{figure}
\begin{picture}(320,30)
\put(0,10){$\varGamma(A_{\infty})=$}
\put(50,10){$\circ$}\put(70,10){$\circ$}\put(90,10){$\circ$}\put(110,10){$\circ$}
\put(54,13){\line(1,0){16}}\put(74,13){\line(1,0){16}}\put(94,13){\line(1,0){16}}\put(114,13){\line(1,0){16}}\put(137,10){$\dotsb$}
\put(50,0){$1$}\put(70,0){$2$}\put(90,0){$3$}\put(110,0){$4$}
\put(170,10){$\varGamma(B_{\infty})=$}
\put(220,10){$\circ$}\put(240,10){$\circ$}\put(260,10){$\circ$}\put(280,10){$\circ$}
\put(224,13){\line(1,0){16}}\put(244,13){\line(1,0){16}}\put(264,13){\line(1,0){16}}\put(284,13){\line(1,0){16}}\put(307,10){$\dotsb$}
\put(220,0){$1$}\put(240,0){$2$}\put(260,0){$3$}\put(280,0){$4$}
\put(230,18){$4$}
\end{picture}\\
\ \\
\begin{picture}(360,40)
\put(0,10){$\varGamma(D_{\infty})=$}
\put(60,0){$\circ$}\put(60,20){$\circ$}\put(80,10){$\circ$}\put(100,10){$\circ$}\put(120,10){$\circ$}
\put(64,3){\line(2,1){16}}\put(64,22){\line(2,-1){16}}
\put(84,13){\line(1,0){16}}\put(104,13){\line(1,0){16}}\put(124,13){\line(1,0){16}}\put(147,10){$\dotsb$}
\put(50,20){$1$}\put(50,0){$2$}\put(80,0){$3$}\put(100,0){$4$}\put(120,0){$5$}
\put(180,10){$\varGamma(E_8)=$}
\put(230,10){$\circ$}\put(250,10){$\circ$}\put(270,10){$\circ$}\put(290,10){$\circ$}\put(310,10){$\circ$}\put(330,10){$\circ$}\put(350,10){$\circ$}\put(270,30){$\circ$}
\put(234,13){\line(1,0){16}}\put(254,13){\line(1,0){16}}\put(274,13){\line(1,0){16}}\put(294,13){\line(1,0){16}}\put(314,13){\line(1,0){16}}\put(334,13){\line(1,0){16}}\put(273,15){\line(0,1){16}}
\put(230,0){$1$}\put(250,0){$3$}\put(270,0){$4$}\put(290,0){$5$}\put(310,0){$6$}\put(330,0){$7$}\put(350,0){$8$}\put(262,30){$2$}
\end{picture}\\
\ \\
\begin{picture}(360,40)
\put(0,10){$\varGamma(F_4)=$}
\put(45,10){$\circ$}\put(65,10){$\circ$}\put(85,10){$\circ$}\put(105,10){$\circ$}
\put(49,13){\line(1,0){16}}\put(69,13){\line(1,0){16}}\put(89,13){\line(1,0){16}}
\put(45,0){$1$}\put(65,0){$2$}\put(85,0){$3$}\put(105,0){$4$}
\put(75,18){$4$}
\put(130,10){$\varGamma(H_4)=$}
\put(180,10){$\circ$}\put(200,10){$\circ$}\put(220,10){$\circ$}\put(240,10){$\circ$}
\put(184,13){\line(1,0){16}}\put(204,13){\line(1,0){16}}\put(224,13){\line(1,0){16}}
\put(180,0){$1$}\put(200,0){$2$}\put(220,0){$3$}\put(240,0){$4$}
\put(190,18){$5$}
\put(265,10){$\varGamma(I_2(m))=$}
\put(325,10){$\circ$}\put(345,10){$\circ$}
\put(329,13){\line(1,0){16}}
\put(325,0){$1$}\put(345,0){$2$}
\put(334,18){$m$}
\end{picture}\\
\ \\
\begin{picture}(320,30)
\put(0,10){$\varGamma(A_{\infty,\infty})=$}
\put(100,10){$\circ$}\put(120,10){$\circ$}\put(140,10){$\circ$}\put(160,10){$\circ$}\put(180,10){$\circ$}\put(200,10){$\circ$}\put(220,10){$\circ$}
\put(65,10){$\dotsb$}\put(100,13){\line(-1,0){16}}\put(104,13){\line(1,0){16}}\put(124,13){\line(1,0){16}}\put(144,13){\line(1,0){16}}\put(164,13){\line(1,0){16}}\put(184,13){\line(1,0){16}}\put(204,13){\line(1,0){16}}\put(224,13){\line(1,0){16}}\put(247,10){$\dotsb$}
\put(95,0){$-3$}\put(115,0){$-2$}\put(135,0){$-1$}\put(160,0){$0$}\put(180,0){$1$}\put(200,0){$2$}\put(220,0){$3$}
\put(270,10){$\supset \varGamma(A_{\infty})$}
\end{picture}
\caption{Some connected Coxeter graphs}
\end{figure}
\end{defn}
By definition, $\varGamma(\mathcal{T}_\infty)$ ($\mathcal{T}=A,B,D$) and $\varGamma(A_{\infty,\infty})$ are Coxeter graphs with countable (infinite) vertex sets.
On the other hand, it is well known that the Coxeter graphs $\varGamma(A_n)$ ($1 \leq n<\infty$), $\varGamma(B_n)$ ($2 \leq n<\infty$), $\varGamma(D_n)$ ($4 \leq n<\infty$), $\varGamma(E_6)$, $\varGamma(E_7)$, $\varGamma(E_8)$, $\varGamma(F_4)$, $\varGamma(H_3)$, $\varGamma(H_4)$ and $\varGamma(I_2(m))$ ($5 \leq m<\infty$) are all the connected Coxeter graphs of finite type (up to isomorphism).
Note that $\varGamma(B_1)=\varGamma(D_1)=\varGamma(A_1)$, while $\varGamma(D_2) \simeq \varGamma(A_1 \times A_1)$ and $\varGamma(D_3) \simeq \varGamma(A_3)$ (but the vertex labels are different).
\subsection{Root systems of Coxeter groups}\label{sec:rootsystem}
For a Coxeter system $(W,S)$, let $\Pi$ be the set of symbols $\alpha_s$ ($s \in S$) and $V$ the vector space over $\mathbb{R}$ containing the set $\Pi$ as a basis.
We define the symmetric bilinear form $\langle \,,\, \rangle$ on $V$ for the basis by
\begin{equation*}
\langle \alpha_s, \alpha_t \rangle=-\cos(\pi/m(s,t)) \text{ if } m(s,t)<\infty,\ \langle \alpha_s, \alpha_t \rangle=-1 \text{ if } m(s,t)=\infty.
\end{equation*}
Then $W$ acts faithfully on the space $V$ by $s \cdot v=v-2\langle \alpha_s,v \rangle \alpha_s$ ($s \in S$, $v \in V$).
Let $\Phi=W \cdot \Pi$, the {\it root system} of $(W,S)$.
The above rule implies that the action of $W$ preserves the bilinear form; as a consequence, any element ({\it root}) of $\Phi$ is a unit vector.
It is a crucial fact that $\Phi$ is a disjoint union of the set $\Phi^+$ of positive roots (i.e.\ roots in which the coefficient of every $\alpha_s \in \Pi$ is $\geq 0$) and the set $\Phi^-=-\Phi^+$ of negative roots.
It is known that the set $\Phi\left[w\right]=\{\gamma \in \Phi^+ \mid w \cdot \gamma \in \Phi^-\}$ characterizes the element $w \in W$; namely,
\begin{equation}
\label{eq:Phiwcharacterizesw}
\text{ if } w,u \in W \text{ and } \Phi\left[w\right]=\Phi\left[u\right], \text{ then } w=u
\end{equation}
(cf.\ Lemma 2.9 of \cite{Nui}, etc.\ for the proof).
Moreover, it is also well known that the cardinality of the set $\Phi\left[w\right]$ is (finite and) equal to the length $\ell(w)$ of $w \in W$ with respect to the generating set $S$.\\ \indent
The {\it reflection} along a root $\gamma=w \cdot \alpha_s \in \Phi$ is defined by $s_\gamma=wsw^{-1} \in W$.
This definition does not depend on the choice of $w$ and $s$, and $s_\gamma$ indeed acts as a reflection on the space $V$; $s_\gamma \cdot v=v-2 \langle \gamma,v \rangle \gamma$ for $v \in V$.
Note that $s_{\alpha_s}=s$ for $s \in S$.
The following fact is easy to show (by the fact that $\Phi=\Phi^+ \sqcup \Phi^-$):
\begin{equation}
\label{eq:reflectionandsimpleroot}
\text{ if } s \in S, \gamma \in \Phi^+ \text{ and } \langle \alpha_s,\gamma \rangle >0, \text{ then } s_\gamma \cdot \alpha_s \in \Phi^-.
\end{equation}
\indent
For $v \in V$, put
\begin{equation*}
v=\sum_{s \in S}(\left[\alpha_s\right]v)\alpha_s \text{ and } {\rm supp}(v)=\{s \in S \mid \left[\alpha_s\right]v \neq 0\}.
\end{equation*}
For $I \subset S$, let $V_I$ be the subspace of $V$ spanned by the set $\Pi_I=\{\alpha_s \mid s \in I\}$ and $\Phi_I=\Phi \cap V_I$ (namely, the set of all $\gamma \in \Phi$ such that ${\rm supp}(\gamma) \subset I$).
Then it is well known that $\Phi_I$ coincides with the root system $W_I \cdot \Pi_I$ of the Coxeter system $(W_I,I)$ (cf.\ Lemma 4 of \cite{Fra-How}, etc.\ for the proof).
This fact yields the following:
\begin{equation}
\label{eq:supportisconnected}
\text{ If } \gamma \in \Phi, \text{ then ($\gamma \in \Phi_{{\rm supp}(\gamma)}$ and so) the set } {\rm supp}(\gamma) \text{ is connected in } \varGamma.
\end{equation}
Moreover, it is well known (cf.\ \cite{Hum}, Section 5.8, Exercise 4, etc.) that:
\begin{equation}
\label{eq:reflectioninparabolic}
\text{ If } I \subset S \text{ and } \gamma \in \Phi, \text{ then } s_\gamma \in W_I \text{ if and only if } \gamma \in \Phi_I.
\end{equation}
For $I \subset S$, let
\begin{equation*}
\begin{split}
I^\perp&=\{s \in S \smallsetminus I \mid st=ts \text{ for all } t \in I\}\\
&=\{s \in S \smallsetminus I \mid s \text{ is adjacent in } \varGamma \text{ to no element of } I\}\\
&=\{s \in S \mid \alpha_s \text{ is orthogonal to every } \alpha_t \in \Pi_I\}.
\end{split}
\end{equation*}
Then we have the following properties:
\begin{align}
\label{eq:rootandsmallerparabolic}
&\text{ If } \gamma \in \Phi^+ \text{ and } {\rm supp}(\gamma) \not\subset I \subset S, \text{ then } w \cdot \gamma \in \Phi^+ \text{ for all } w \in W_I.\\
\label{eq:rootandadjacentelement}
&\text{ If } \gamma \in \Phi, I={\rm supp}(\gamma) \text{ and } s \in S \smallsetminus (I \cup I^\perp), \text{ then } {\rm supp}(s \cdot \gamma)=I \cup \{s\}.
\end{align}
(For (\ref{eq:rootandsmallerparabolic}), take some $t \in {\rm supp}(\gamma) \smallsetminus I$, then $w \cdot \gamma$ has the same (positive) coefficient of $\alpha_t$ as $\gamma$.
For (\ref{eq:rootandadjacentelement}), note that $\langle \alpha_s,\gamma \rangle<0$ by the hypothesis.)\\ \indent
For $I \subset S$ and $w \in W$, let $\Phi_I^+=\Phi_I \cap \Phi^+$, $\Phi_I^-=\Phi_I \cap \Phi^-$ and $\Phi_I\left[w\right]=\Phi_I \cap \Phi\left[w\right]$.
\begin{lem}
\label{lem:lemmaforreflectiondecomposition}
Let $w \in W$, $I,J \subset S$ and suppose that $I \cap J=\emptyset$, $w \cdot \Pi_I=\Pi_I$ and $w \cdot \Pi_J \subset \Phi^-$.
Then $\Phi_{I \cup J}\left[w\right]=\Phi_{I \cup J}^+ \smallsetminus \Phi_I$.
\end{lem}
\begin{proof}
Let $\gamma \in \Phi_{I \cup J}^+$ such that $\left[\alpha_s\right]\gamma>0$ for at least one $s \in J$ (note that $w \cdot \alpha_s \in \Phi^-$).
Now if $w \cdot \alpha_s \in \Phi_I^-$, then $\alpha_s=w^{-1} \cdot (w \cdot \alpha_s)$ must be a linear combination of $\Pi_I$ (since $w \cdot \Pi_I=\Pi_I$), but this is impossible.
Thus we have $\left[\alpha_t\right](w \cdot \alpha_s)<0$ for some $t \in S \smallsetminus I$.
Moreover, the hypothesis implies that $\left[\alpha_t\right](w \cdot \alpha_{s'})=0$ for all $s' \in I$ and $\left[\alpha_t\right](w \cdot \alpha_{s'}) \leq 0$ for all $s' \in J$.
Thus we have
\begin{equation*}
\begin{split}
\left[\alpha_t\right](w \cdot \gamma)&=\left[\alpha_t\right]\biggl(w \cdot \sum_{s' \in I \cup J}(\left[\alpha_{s'}\right]\gamma)\alpha_{s'}\biggr)\\
&=\sum_{s' \in I \cup J}(\left[\alpha_{s'}\right]\gamma)\left[\alpha_t\right](w \cdot \alpha_{s'}) \leq \left[\alpha_s\right]\gamma\left[\alpha_t\right](w \cdot \alpha_s)<0.
\end{split}
\end{equation*}
Hence the claim holds, since $w \cdot \Phi_I^+ \subset \Phi^+$ by the hypothesis.
\end{proof}
\begin{defn}
\label{defn:oddCoxetergraph}
For a Coxeter system $(W,S)$, we define the {\it odd Coxeter graph} $\varGamma^{\rm odd}$ of $(W,S)$ to be the subgraph of $\varGamma$ obtained by removing all edges labelled by an even number or $\infty$.
\end{defn}
It is well known (cf.\ \cite{Hum}, Section 5.3, Exercise, etc.) that, for $s,t \in S$,
\begin{equation}
\label{eq:oddCoxetersubgraph}
\alpha_t \in W \cdot \alpha_s \text{ if and only if } s,t \text{ are in the same connected component of $\varGamma^{\rm odd}$}.
\end{equation}
Moreover, the following lemma is deduced immediately from the definition that all fundamental relations of $W$ are of the form $(st)^{m(s,t)}=1$ ($s,t \in S$).
\begin{lem}
\label{lem:charofWtopm1}
Any $f \in {\rm Hom}(W,\{\pm 1\})$ assigns the same value to every vertex $s \in S$ of a connected component of $\varGamma^{\rm odd}$.
Conversely, any mapping $S \to \{\pm 1\}$ having this property extends uniquely to a homomorphism $W \to \{\pm 1\}$.
\end{lem}
\subsection{Reflection decompositions of longest elements}\label{sec:decompositionofw_0}
If $W_I$ is a finite parabolic subgroup of a Coxeter group $W$, then let $w_0(I)$ denote the {\it longest element} of $W_I$.
This element is an involution and maps the set $\Pi_I$ onto $-\Pi_I$, so that there is an involutive graph automorphism $\sigma_I$ of the Coxeter graph $\varGamma_I$ such that
\begin{equation*}
w_0(I) \cdot \alpha_s=-\alpha_{\sigma_I(s)} \text{ for all } s \in I.
\end{equation*}
It is well known that, for an irreducible Coxeter system $(W,S)$, we have $Z(W) \neq 1$ if and only if $W \simeq W(\mathcal{T})$ for one of $\mathcal{T}=A_1$, $B_n$ ($n<\infty$), $D_k$ ($k \geq 4$ even), $E_7$, $E_8$, $F_4$, $H_3$, $H_4$ and $I_2(m)$ ($m \geq 6$ even).
This condition is also equivalent to that $|W|<\infty$ and $\sigma_{S}={\rm id}_S$.
Moreover, $Z(W)=\{1,w_0(S)\}$ if $Z(W) \neq 1$, while $\sigma_S$ is determined as the unique non-identical automorphism of $\varGamma$ whenever $W$ is finite, irreducible and $Z(W)=1$.
Note that any automorphism $\tau \in {\rm Aut}(\varGamma)$ induces naturally an automorphism of $W$, which maps each element $w_0(I)$ to $w_0(\tau(I))$.\\ \indent
In the paper \cite{Deo}, Deodhar established a method (in the proof of Theorem 5.4) for decomposing any involution $w \in W$ as a product of commuting reflections.
From now, we apply this method and then obtain a decomposition of any longest element $w_0(I)$, which we call here a {\it reflection decomposition}.
First, to each finite irreducible Coxeter system $(W,S)=(W(\mathcal{T}),S(\mathcal{T}))$ of type $\mathcal{T}$, we associate a (or two) positive root(s) $\widetilde{\alpha}_{\mathcal{T}}=\widetilde{\alpha}_{\mathcal{T}}^{(1)}$ (and $\widetilde{\alpha}_{\mathcal{T}}^{(2)}$), as follows (where we abbreviate $c_1\alpha_1+c_2\alpha_2+ \dotsb +c_n\alpha_n \in V$ to $(c_1,c_2,\dotsc ,c_n)$ in some cases):
\begin{equation*}
\begin{split}
&\widetilde{\alpha}_{A_n}=\sum_{i=1}^{n}\alpha_i\quad (1 \leq n<\infty),\ \widetilde{\alpha}_{D_n}=\alpha_1+\alpha_2+\sum_{i=3}^{n-1}2\alpha_i+\alpha_n\quad (4 \leq n<\infty),\\
&\widetilde{\alpha}_{B_n}^{(1)}=\alpha_1+\sum_{i=2}^{n}\sqrt{\,2\,}\alpha_i,\ \widetilde{\alpha}_{B_n}^{(2)}=\sqrt{\,2\,}\alpha_1+\sum_{i=2}^{n-1}2\alpha_i+\alpha_n\quad (2 \leq n<\infty),\\
&\widetilde{\alpha}_{E_6}=(1,2,2,3,2,1),\ \widetilde{\alpha}_{E_7}=(2,2,3,4,3,2,1),\ \widetilde{\alpha}_{E_8}=(2,3,4,6,5,4,3,2),\\
&\widetilde{\alpha}_{F_4}^{(1)}=(2,3,2\sqrt{\,2\,},\sqrt{\,2\,}),\ \widetilde{\alpha}_{F_4}^{(2)}=(\sqrt{\,2\,},2\sqrt{\,2\,},3,2),\\
&\widetilde{\alpha}_{H_3}=(c+1,2c,c),\ \widetilde{\alpha}_{H_4}=(3c+2,4c+2,3c+1,2c)\quad (\text{where } c=2\cos\frac{\pi}{5}),\\
&\widetilde{\alpha}_{I_2(m)}=\frac{1}{2\sin(\pi/2m)}\alpha_1+\frac{1}{2\sin(\pi/2m)}\alpha_2\quad (m \geq 5 \text{ odd}),\\
&\widetilde{\alpha}_{I_2(m)}^{(i)}=\frac{\cos(\pi/m)}{\sin(\pi/m)}\alpha_i+\frac{1}{\sin(\pi/m)}\alpha_{3-i}\quad (m \geq 5 \text{ even}, i=1,2).
\end{split}
\end{equation*}
To check that each of these is actually a root of $(W(\mathcal{T}),S(\mathcal{T}))$, note the equality $c^2=c+1$ and the following formula for the root system of type $I_2(m)$:
\begin{equation*}
\begin{split}
&\text{ If } w=(\dotsm s_2s_1s_2) \in W(I_2(m)) \text{ ($k$ elements), then}\\
&w \cdot \alpha_1=
\begin{cases}
\frac{\sin(k\pi/m)}{\sin(\pi/m)}\alpha_1+\frac{\sin((k+1)\pi/m)}{\sin(\pi/m)}\alpha_2 & \text{ if } k \text{ is odd},\medskip\\
\frac{\sin((k+1)\pi/m)}{\sin(\pi/m)}\alpha_1+\frac{\sin(k\pi/m)}{\sin(\pi/m)}\alpha_2 & \text{ if } k \text{ is even}.
\end{cases}
\end{split}
\end{equation*}
\indent
For example, we have
\begin{equation*}
\begin{split}
&\widetilde{\alpha}_{F_4}^{(1)}=s_1s_2s_3s_4s_2s_3s_2 \cdot \alpha_1,\quad \widetilde{\alpha}_{F_4}^{(2)}=s_4s_3s_2s_1s_3s_2s_3 \cdot \alpha_4,\\
&\widetilde{\alpha}_{H_3}=s_2s_1s_2s_1s_3s_2 \cdot \alpha_1,\quad \widetilde{\alpha}_{H_4}=s_4s_3s_2s_1s_2s_1s_3s_2s_1s_4s_3s_2s_1s_2s_3s_4 \cdot \widetilde{\alpha}_{H_3},\\
&\widetilde{\alpha}_{I_2(2k+1)}=(\dotsm s_2s_1s_2) \cdot \alpha_1\ (k \text{ elements}),\ \widetilde{\alpha}_{I_2(4k)}^{(i)}=(s_{3-i}s_i)^{k-1}s_{3-i} \cdot \alpha_i.
\end{split}
\end{equation*}
By (\ref{eq:oddCoxetersubgraph}), if $\mathcal{T} \neq B_n$, $F_4$, $I_2(m)$ ($m$ even), then $\Phi$ consists of a single orbit $W(\mathcal{T}) \cdot \alpha_1$ (and so it contains $\widetilde{\alpha}_{\mathcal{T}}$).
On the other hand, if $\mathcal{T}=B_n$, $F_4$ or $I_2(4k)$, then (\ref{eq:oddCoxetersubgraph}) implies that $\Phi$ consists of two orbits (namely, $W \cdot \alpha_1$ and $W \cdot \alpha_2$ if $\mathcal{T}=B_n$, $I_2(4k)$, and $W \cdot \alpha_1$ and $W \cdot \alpha_4$ if $\mathcal{T}=F_4$).
In these case, $\widetilde{\alpha}_{\mathcal{T}}^{(1)}$ lies in the orbit $W \cdot \alpha_1$ and $\widetilde{\alpha}_{\mathcal{T}}^{(2)}$ lies in the other one.\\ \indent
In contrast with the above cases, if $\mathcal{T}=I_2(4k+2)$, then $\Phi$ consists of two orbits $W(\mathcal{T}) \cdot \alpha_1$ and $W(\mathcal{T}) \cdot \alpha_2$, and now we have $\widetilde{\alpha}_{\mathcal{T}}^{(1)} \in W(\mathcal{T}) \cdot \alpha_2$ (and $\widetilde{\alpha}_{\mathcal{T}}^{(2)}$ lies in the other orbit).
In fact, we have $\widetilde{\alpha}_{I_2(4k+2)}^{(i)}=(s_{3-i}s_i)^k \cdot \alpha_{3-i}$ for $i=1,2$.\\ \indent
To simplify the description, we denote the reflection along the root $\widetilde{\alpha}_{\mathcal{T}}^{(i)}$ by $\widetilde{r}(\mathcal{T},i)$.
If we have only one root $\widetilde{\alpha}_{\mathcal{T}}^{(i)}$, namely $\mathcal{T} \neq B_n$, $F_4$, $I_2(m)$ ($m$ even), then we also write $\widetilde{r}(\mathcal{T})=\widetilde{r}(\mathcal{T},1)$.
\begin{rem}
\label{rem:highestroots}
By the above observation, if $\mathcal{T}=B_n$, $F_4$ or $I_2(4k)$, then $\widetilde{r}(\mathcal{T},1)$ is conjugate to $s_1$, and $\widetilde{r}(\mathcal{T},2)$ is conjugate to $s_2$ (if $\mathcal{T}=B_n$ or $I_2(4k)$) or to $s_4$ (if $\mathcal{T}=F_4$).
On the other hand, if $\mathcal{T}=I_2(4k+2)$, then $\widetilde{r}(\mathcal{T},1)$, $\widetilde{r}(\mathcal{T},2)$ are conjugate to $s_2$, $s_1$, respectively.
\end{rem}
\begin{lem}
\label{lem:highestroots}
{\bf (i)} If $\mathcal{T} \neq A_n$ ($n \geq 2$), $I_2(m)$ ($m$ odd), then for the root $\widetilde{\alpha}_{\mathcal{T}}^{(i)}$, there is an index $N(\mathcal{T},i)$ such that $\langle \widetilde{\alpha}_{\mathcal{T}}^{(i)},\alpha_j \rangle=0$ for all $j \neq N(\mathcal{T},i)$.
Moreover, we have $\langle \widetilde{\alpha}_{\mathcal{T}}^{(i)},\alpha_{N(\mathcal{T},i)} \rangle >0$ and $\Phi\left[\!\right.\widetilde{r}(\mathcal{T},i)\left.\!\right]=\Phi^+ \smallsetminus \Phi_{S(\mathcal{T}) \smallsetminus \{s_{N(\mathcal{T},i)}\}}$.
(If we have only one root $\widetilde{\alpha}_{\mathcal{T}}^{(i)}$, then we also write $N(\mathcal{T})=N(\mathcal{T},1)$.)\\
{\bf (ii)} If $\mathcal{T}=A_n$ ($n \geq 2$) or $I_2(m)$ ($m$ odd), then there are two indices $N_1(\mathcal{T}),N_2(\mathcal{T})$ such that $\langle \widetilde{\alpha}_{\mathcal{T}},\alpha_{N_j(\mathcal{T})} \rangle >0$ for $j=1,2$ and $\langle \widetilde{\alpha}_{\mathcal{T}},\alpha_j \rangle=0$ for all $j \neq N_1(\mathcal{T}),N_2(\mathcal{T})$.
Moreover, we have $\Phi\left[\!\right.\widetilde{r}(\mathcal{T},i)\left.\!\right]=\Phi^+ \smallsetminus \Phi_{S(\mathcal{T}) \smallsetminus \{s_{N_1(\mathcal{T})},s_{N_2(\mathcal{T})}\}}$.
\end{lem}
\begin{proof}
{\bf (i)} The first claim follows from a direct computation, by putting
\begin{equation*}
\begin{split}
&N(A_1)=1,\quad N(B_n,1)=n,\quad N(B_n,2)=n-1,\quad N(D_n)=n-1,\\
&N(E_6)=2,\quad N(E_7)=1,\quad N(E_8)=8,\quad N(F_4,1)=1,\quad N(F_4,2)=4,\\
&N(H_3)=2,\quad N(H_4)=4,\quad N(I_2(2k),1)=2,\quad N(I_2(2k),2)=1.
\end{split}
\end{equation*}
For the second one, expand the equality $\langle \widetilde{\alpha}_{\mathcal{T}}^{(i)},\widetilde{\alpha}_{\mathcal{T}}^{(i)} \rangle=1$ and use the first claim.
Now the third one follows from (\ref{eq:reflectionandsimpleroot}) and Lemma \ref{lem:lemmaforreflectiondecomposition}.\\
{\bf (ii)} The former claim also follows from a direct computation, by putting
\begin{equation*}
N_1(A_n)=1,\ N_2(A_n)=n,\quad N_1(I_2(2k+1))=1,\ N_2(I_2(2k+1))=2.
\end{equation*}
The remaining proof is similar to (i).
\end{proof}
Now Deodhar's method can be described, for the element $w_0(I)$, as follows:\medskip\\ \indent
{\bf (I)} If $I=\emptyset$, then this algorithm finishes with the (trivial) decomposition $w_0(I)=1$.
If $I \neq \emptyset$, choose an irreducible component $J$ of $I$.
Let $J=S(\mathcal{T})$.\\ \indent
{\bf (II)} If $\mathcal{T} \neq A_n$ ($n \geq 2$), $I_2(m)$ ($m$ odd), take the (or one of the two) root(s) $\widetilde{\alpha}_{\mathcal{T}}^{(i)}$.
By Lemma \ref{lem:highestroots} (i), $\widetilde{r}(\mathcal{T},i)$ commutes with all elements of $K=I \smallsetminus \{s_{N(\mathcal{T},i)}\}$, and we have $w_0(I)=\widetilde{r}(\mathcal{T},i)w_0(K)$ (since both sides map $\Pi_I$ into $\Phi^-$; cf.\ (\ref{eq:Phiwcharacterizesw})).
Then apply this algorithm inductively to the (smaller) set $K$.\\ \indent
{\bf (III)} If $\mathcal{T}=A_n$ ($n \geq 2$) or $I_2(m)$ ($m$ odd), then similarly, $\widetilde{r}(\mathcal{T})$ commutes with all elements of $K=I \smallsetminus \{s_{N_1(\mathcal{T})},s_{N_2(\mathcal{T})}\}$ and $w_0(I)=\widetilde{r}(\mathcal{T},1)w_0(K)$ by Lemma \ref{lem:highestroots} (ii).
Then apply this algorithm inductively to the (smaller) set $K$.\medskip\\ \indent
By collecting the subset $K \subset I$ appearing in the step (II) or (III) of every turn, we obtain a decreasing sequence ($K_0=I$,) $K_1,\dotsc ,K_{r-1},K_r=\emptyset$.
We call this a {\it generator sequence} (of length $r$) for the set $I$.
\begin{exmp}
Let $(W,S)=(W(D_n),S(D_n))$.
By using a reflection decomposition of $w_0(S(D_i))$, we compute the root $s_{i+1}w_0(S(D_i))s_{i+1} \cdot \alpha_i$ ($3 \leq i<n$).
First, assume that $i$ is odd.
By the algorithm, we have a decomposition
\begin{equation*}
w_0(S(D_i))=\widetilde{r}(D_i)s_i\widetilde{r}(D_{i-2})s_{i-2} \dotsm \widetilde{r}(D_5)s_5\widetilde{r}(D_3)s_3
\end{equation*}
(where we put $\widetilde{r}(D_3)=s_{\alpha_1+\alpha_2+\alpha_3}$; note that $\varGamma(D_3) \simeq \varGamma(A_3)$).
The corresponding generator sequence is
\begin{equation*}
\begin{split}
S(D_{i-2}) \cup \{s_i\},\ S(D_{i-2}),&\ S(D_{i-4}) \cup \{s_{i-2}\},\ S(D_{i-4}),\dotsc\\
&\dotsc ,\ S(D_5),\ S(D_3) \cup \{s_5\},\ S(D_3),\ \{s_3\},\ \emptyset.
\end{split}
\end{equation*}
Now since $\sigma_{S(D_i)}(s_i)=s_i$, we have
\begin{equation*}
\begin{split}
w_0(S(D_i))s_{i+1}& \cdot \alpha_i=w_0(S(D_i)) \cdot (\alpha_i+\alpha_{i+1})\\
&=w_0(S(D_i)) \cdot \alpha_i+w_0(S(D_i)) \cdot \alpha_{i+1}=-\alpha_i+w_0(S(D_i)) \cdot \alpha_{i+1}.
\end{split}
\end{equation*}
Since all the reflections except $\widetilde{r}(D_i)$, $s_i$ in the decomposition fix the root $\alpha_{i+1}$, and all roots corresponding to the reflections are orthogonal (by definition), we have
\begin{equation*}
w_0(S(D_i)) \cdot \alpha_{i+1}=\alpha_{i+1}-2\langle \widetilde{\alpha}_{D_i}, \alpha_{i+1}\rangle \widetilde{\alpha}_{D_i}-2\langle \alpha_i, \alpha_{i+1}\rangle \alpha_i=\widetilde{\alpha}_{D_{i+1}}
\end{equation*}
(where we put $\widetilde{\alpha}_{D_3}=\alpha_1+\alpha_2+\alpha_3$).
Thus we have
\begin{equation*}
s_{i+1}w_0(S(D_i))s_{i+1} \cdot \alpha_i=s_{i+1} \cdot (\widetilde{\alpha}_{D_{i+1}}-\alpha_i)=\widetilde{\alpha}_{D_i}.
\end{equation*}
\indent
On the other hand, if $i \geq 3$ is even, then we have a different decomposition
\begin{equation*}
w_0(S(D_i))=\widetilde{r}(D_i)s_i\widetilde{r}(D_{i-2})s_{i-2} \dotsm \widetilde{r}(D_4)s_4s_2s_1.
\end{equation*}
However, we obtain the same result; namely, we have
\begin{equation*}
s_{i+1}w_0(S(D_i))s_{i+1} \cdot \alpha_i=\widetilde{\alpha}_{D_i}.
\end{equation*}
\end{exmp}
By a similar argument, it can be checked that $s_{i+1}w_0(S(D_i))s_{i+1}$ ($i \geq 3$) maps the roots $\alpha_{i+1}$, $\alpha_i$, $\alpha_j$ ($j<i$) to $-\widetilde{\alpha}_{D_{i+1}}$, $\widetilde{\alpha}_{D_i}$, $-\alpha_{j'}$ (where $j'$ is the index such that $s_{j'}=\sigma_{S(D_i)}(s_j)$) respectively.
The element $w_0(S(D_{i-1}))w_0(S(D_i))w_0(S(D_{i+1}))$ has the same property.
Thus we have
\begin{equation*}
s_{i+1}w_0(S(D_i))s_{i+1}=w_0(S(D_{i-1}))w_0(S(D_i))w_0(S(D_{i+1}))\quad (i \geq 3).
\end{equation*}
Similarly, we have the following relations:
\begin{equation*}
\begin{split}
s_{i+1}w_0(S(B_i))s_{i+1}&=w_0(S(B_{i-1}))w_0(S(B_i))w_0(S(B_{i+1}))\quad (i \geq 2),\\
s_2s_1s_2&=s_1w_0(S(B_2)),\\
s_3w_0(S(D_2))s_3&=w_0(S(D_2))w_0(S(D_3)),\\
w_0(S(D_i))w_0(S(D_j))&=w_0(S(D_j))w_0(S(D_i))\quad (2 \leq i<j),\\
s_1w_0(S(D_{2k+1}))s_1&=s_2w_0(S(D_{2k+1}))s_2=w_0(S(D_2))w_0(S(D_{2k+1})).
\end{split}
\end{equation*}
(The last row follows from the relations $w_0(S(D_{2k+1})) \cdot \alpha_i=-\alpha_{3-i}$ ($i=1,2$).)
Moreover, note that $w_0(S(B_i)) \in Z(W(B_i))$ and $w_0(S(D_{2k})) \in Z(W(D_{2k}))$, and $w_0(S(D_{2k+1}))$ commutes with all $s_j$ ($3 \leq j \leq 2k+1$).\\ \indent
By these relations, we have the following:
\begin{lem}
\label{lem:defofpropersubgroupG}
(See Definition \ref{defn:Coxetergraphs} for notations.)\\
{\bf (i)} Let $1 \leq n \leq \infty$.
Then the subgroup $G_{B_n}$ of $W(B_n)$ generated by all $w_0(S(B_i))$ ($1 \leq i \leq n$, $i<\infty$) is normal in $W(B_n)$.\\
{\bf (ii)} Let $1 \leq n \leq \infty$.
Then the smallest normal subgroup $G_{D_n}$ of $W(D_n)$ containing all $w_0(S(D_{2k}))$ ($1 \leq k<\infty$, $2k \leq n$) is the subgroup generated by all $w_0(S(D_i))$ ($2 \leq  i \leq n$, $i<\infty$).\\
{\bf (iii)} Moreover, each of the above normal subgroups is an elementary abelian $2$-group with the generating set given there as the basis.
\end{lem}
These normal subgroups $G_{B_n}$, $G_{D_n}$ will appear in later sections.
\subsection{Direct product decompositions of finite Coxeter groups}\label{sec:decompositionoffiniteW}
Owing to the reflection decomposition given in Section \ref{sec:decompositionofw_0}, we can determine easily which finite irreducible Coxeter groups have the center as a nontrivial direct factor.
(This is never a new result, but we restate it here since the result is used in later sections.)\\ \indent
For a Coxeter system $(W,S)$, let $W^+$ denote the normal subgroup of $W$ (of index two) consisting of elements of even length.
This coincides with the kernel of the map ${\rm sgn} \in {\rm Hom}(W,\{\pm 1\})$ such that ${\rm sgn}(w)=(-1)^{\ell(w)}$.
Since any reflection in $W$ has odd length, the following lemma follows from (the proof of) Lemma \ref{lem:conditionforvanishingcenter}:
\begin{lem}
\label{lem:lemmafordirectfactorZ}
If $(W,S)$ is a finite irreducible Coxeter system and $Z(W) \neq 1$, then we have $W=Z(W) \times W^+$ if and only if some (or equivalently, any) generator sequence for $S$ (cf.\ Section \ref{sec:decompositionofw_0}) has odd length.
\end{lem}
\begin{thm}
\label{thm:Zisdirectfactor}
Let $(W,S)$ be an irreducible Coxeter system such that $Z(W) \neq 1$ (so that $|W|<\infty$).
Then $Z(W)$ ($\simeq W(A_1)$) is a proper direct factor of $W$ if and only if $W \simeq W(\mathcal{T})$ for $\mathcal{T}=B_{2k+1}$, $I_2(4k+2)$ ($k \geq 1$), $E_7$ or $H_3$.
In the first two cases, $W$ is isomorphic to $W(A_1) \times W(D_{2k+1})$, $W(A_1) \times W(I_2(2k+1))$ respectively.
In the last two cases, we have $W=Z(W) \times W^+$.
\end{thm}
\begin{proof}
Note that $Z(W) \simeq \{\pm 1\}$ by the hypothesis.
Since $Z(W(A_1))=W(A_1)$, we may assume $W \neq W(A_1)$.\\ \indent
{\bf Case 1. $W=W(B_n)$ ($n \geq 2$):}\ 
First, we have ${\rm Hom}(W,\{\pm 1\})=\{1,{\rm sgn},\varepsilon_1,\varepsilon_2\}$ by Lemma \ref{lem:charofWtopm1}, where $1$ denotes the trivial map, $\varepsilon_1(s_1)=-1$, $\varepsilon_1(s_i)=1$, $\varepsilon_2(s_1)=1$ and $\varepsilon_2(s_i)=-1$ ($i \neq 1$).
Now we consider the following reflection decomposition:
\begin{equation*}
w_0(S)=\widetilde{r}(B_n,1)\widetilde{r}(B_{n-1},1) \dotsm \widetilde{r}(B_2,1)s_1.
\end{equation*}
By Remark \ref{rem:highestroots}, each reflection $\widetilde{r}(B_k,1)$ is conjugate to $s_1$.
This implies that any expression of $\widetilde{r}(B_k,1)$ as a product of generators contains an odd number of $s_1$ and an even number of $s_i$ ($i \neq 1$).
Thus we have
\begin{equation*}
{\rm sgn}(\widetilde{r}(B_k,1))=\varepsilon_1(\widetilde{r}(B_k,1))=-1 \text{ and } \varepsilon_2(\widetilde{r}(B_k,1))=1.
\end{equation*}
If $n$ is even, then all $f \in {\rm Hom}(W,\{\pm 1\})$ maps $w_0(S)$ to $1$ by the above property.
Thus by Lemma \ref{lem:conditionforvanishingcenter}, $Z(W)$ is not a direct factor.\\ \indent
On the other hand, if $n$ is odd, then we have $\varepsilon_1(w_0(S))=-1$ and so $W=Z(W) \times \ker \varepsilon_1$ by the proof of Lemma \ref{lem:conditionforvanishingcenter}.
Note that $\ker \varepsilon_1$ consists of elements in which $s_1$ appears an even number of times.
Since $s_1$ commutes with all $s_i$ ($3 \leq i \leq n$), it can be deduced directly that $\ker \varepsilon_1$ is generated by $s'_1=s_1s_2s_1$ and all $s'_i=s_i$ ($2 \leq i \leq n$).
Moreover, $\ker \varepsilon_1$ forms a Coxeter group of type $D_n$; in fact, $s'_1,\dotsc ,s'_n$ satisfy the fundamental relations of type $D_n$ (so that $\ker \varepsilon_1$ is a quotient of $W(D_n)$), while the order $|W(B_n)|/2$ of $\ker \varepsilon_1$ coincides with $|W(D_n)|$.
Hence the claim holds in this case.\\ \indent
{\bf Case 2. $W=W(\mathcal{T})$ for $\mathcal{T}=D_{2k}$ ($k \geq 2$), $E_7$, $E_8$, $H_3$, $H_4$:}\ 
Since $\varGamma^{\rm odd}$ is connected in this case, we have ${\rm Hom}(W,\{\pm 1\})=\{1,{\rm sgn}\}$ by Lemma \ref{lem:charofWtopm1}.
Thus the claim follows from Lemmas \ref{lem:conditionforvanishingcenter} and \ref{lem:lemmafordirectfactorZ}, by taking the following generator sequence for $S$ (where we abbreviate the set $\{s_{i_1},s_{i_2},\dotsc ,s_{i_r}\}$ to $i_1i_2 \dotsm i_r$):
\begin{equation*}
\begin{cases}
S(D_{2k-2}) \cup \{s_{2k}\},\ S(D_{2k-2}),\dotsc ,\ S(D_4),\ 124,\ 12,\ 1,\ \emptyset & \text{ if } \mathcal{T}=D_{2k},\\
S(E_7),\ 234567,\ 23457,\ 2345,\ 235,\ 23,\ 2,\ \emptyset & \text{ if } \mathcal{T}=E_8,\\
234567,\ 23457,\ 2345,\ 235,\ 23,\ 2,\ \emptyset & \text{ if } \mathcal{T}=E_7,\\
S(H_3),\ 13,\ 1,\ \emptyset & \text{ if } \mathcal{T}=H_4,\\
13,\ 1,\ \emptyset & \text{ if } \mathcal{T}=H_3
\end{cases}
\end{equation*}
(note that the first sequence consists of $2k$ terms).\\ \indent
{\bf Case 3. $W=W(F_4)$:}\ 
We have a generator sequence $234$, $23$, $2$, $\emptyset$ for $S$ and the corresponding decomposition of $w_0(S)$ into four reflections, all of which are conjugate to $s_1$ and $s_2$ (cf.\ Remark \ref{rem:highestroots}).
This (and Lemma \ref{lem:charofWtopm1}) implies that any $f \in {\rm Hom}(W,\{\pm 1\})$ maps all the four reflections to the same element $f(s_1)$, so that $f(w_0(S))=1$.
Hence the claim follows from Lemma \ref{lem:conditionforvanishingcenter}.\\ \indent
{\bf Case 4. $W=W(I_2(2k))$ ($k \geq 3$):}\ 
We have a reflection decomposition $w_0(S)=\widetilde{r}(I_2(2k),1)s_1$.
If $k$ is even, then $\widetilde{r}(I_2(2k),1)$ is conjugate to $s_1$ (cf.\ Remark \ref{rem:highestroots}).
Now by a similar argument to the previous case, any $f \in {\rm Hom}(W,\{\pm 1\})$ maps $w_0(S)$ to $1$.
Thus $Z(W)$ is not a direct factor by Lemma \ref{lem:conditionforvanishingcenter}.\\ \indent
On the other hand, if $k$ is odd, then $\widetilde{r}(I_2(2k),1)$ is conjugate to $s_2$ (cf.\ Remark \ref{rem:highestroots}).
Thus $\varepsilon_1 \in {\rm Hom}(W,\{\pm 1\})$ ($\varepsilon(s_1)=-1$, $\varepsilon(s_2)=1$) sends $w_0(S)$ to $-1$, so that $W=Z(W) \times \ker \varepsilon_1$ by the proof of Lemma \ref{lem:conditionforvanishingcenter}.
Moreover, $\ker \varepsilon_1$ is generated by two reflections $s_1s_2s_1$ and $s_2$, and so $\ker \varepsilon_1$ is a Coxeter system of type $I_2(k)$ (since $s_1s_2s_1s_2$ has order $k$).
Hence the claim holds in all cases.
\end{proof}
Since the groups $W(E_7)^+$ and $W(H_3)^+$ are known to (be isomorphic to) the well-examined simple groups $S_6(2)$ and $A_5$ respectively (cf.\ \cite{Hum}, Sections 2.12--13, etc.), we omit the proof of the following properties of these groups.
Note that these properties can also be proved by using Theorems \ref{thm:Zisdirectfactor} and \ref{thm:indecomposabilityofW} below.
\begin{lem}
\label{lem:evensubgroup}
Let $G=W(\mathcal{T})^+$, $\mathcal{T} \in \{E_7,H_3\}$.
Then $G$ has trivial center, is directly indecomposable and is generated by involutions.
Moreover, $G$ is not isomorphic to a Coxeter group.
\end{lem}
%
\subsection{Notes on normalizers in Coxeter groups}\label{sec:noteonnormalizer}
In this subsection, we summarize some properties of normalizers $N_W(W_I)$ of parabolic subgroups $W_I$ in Coxeter groups $W$.
In the paper \cite{Bri-How} (or \cite{How}, for the case $|W|<\infty$), the structure of $N_W(W_I)$ is well examined so that we can in fact determine the precise structure of the normalizer.
In particular, here we use the following results in those papers:
\begin{prop}[\cite{Bri-How}, Proposition 2.1]
\label{prop:Brink-Howlett_decomposition}
If $I \subset S$, then $N_W(W_I)$ is the semidirect product of $W_I$ by the group $N_I=\{w \in W \mid w \cdot \Pi_I=\Pi_I\}$.
\end{prop}
\begin{prop}[\cite{Bri-How}, remarked between Theorems A and B]
\label{prop:Nofinfiniteirreducible}
If $I \subset J \subset S$ and $W_I$ is an infinite irreducible component of $W_J$, then $N_W(W_J) \subset W_{I \cup I^\perp}$.
\end{prop}
By using these, we can prove the following corollary.
(This is also a consequence of a result in \cite{Bri-How}, but we include the proof here since it is sufficiently short.)\\ \indent
In the proof, we also use the following result.
(This result was originally given by Deodhar \cite{Deo}, in the proof of Proposition 4.2, for the case $|S|<\infty$.
See also \cite{Nui}, Proposition 2.14, etc. for the case $|S|=\infty$.)
\begin{prop}
\label{prop:Deodhar_infiniteroots}
If $(W,S)$ is irreducible and $|W|=\infty$, then $|\Phi \smallsetminus \Phi_I|=\infty$ for all proper subsets $I \subset S$.
\end{prop}
\begin{cor}
\label{cor:Nofmaximalparabolic}
Let $s \in S$ and $I=S \smallsetminus \{s\}$.\\
{\bf (i)} If $1 \neq w \in N_I$, then $\Phi\left[w\right]=\Phi^+ \smallsetminus \Phi_I$.
Hence by (\ref{eq:Phiwcharacterizesw}), such an element $w$ is unique if it exists.\\
{\bf (ii)} If $|W|<\infty$ and $w_0(S) \in N_W(W_I)$, then $N_W(W_I)=W_I \rtimes \{1,w_0(S)\}$.\\
{\bf (iii)} If $(W,S)$ is irreducible and $|W|=\infty$, then $N_I=1$ and $N_W(W_I)=W_I$.
\end{cor}
\begin{proof}
{\bf (i)} In this case, we have $w \cdot \alpha_s \in \Phi^-$ (otherwise, we have $w \cdot \Phi^+ \subset \Phi^+$ but this is a contradiction).
Now the claim follows from Lemma \ref{lem:lemmaforreflectiondecomposition}.\\
{\bf (ii)} Note that $w_0(S) \not\in W_I$, while $|N_I| \leq 2$ by (i).
Thus by Proposion \ref{prop:Brink-Howlett_decomposition}, $N_W(W_I)$ is generated by $W_I$ and $w_0(S)$.
Now the claim holds, since $w_0(S)^2=1$.\\
{\bf (iii)} In this case, we have $|\Phi^+ \smallsetminus \Phi_I|=\infty$ by Proposition \ref{prop:Deodhar_infiniteroots}.
Thus we have $N_I=1$ by (i), since the set $\Phi\left[w\right]$ is always finite.
Hence the claim holds.
\end{proof}
Owing to this description, we have the following:
\begin{cor}
\label{cor:Nforspecialcase}
{\bf (i)} If $W=W(B_n)$, $2 \leq n<\infty$, then $\bigcap_{i=1}^{n-1}N_W(W_{S(B_i)})=G_{B_n}$.\\
{\bf (ii)} If $W=W(D_n)$, $3 \leq n<\infty$, then $\bigcap_{i=2}^{n-1}N_W(W_{S(D_i)})=G_{D_n} \rtimes \langle s_1 \rangle$.
\end{cor}
\begin{proof}
Note that, by Lemma \ref{lem:defofpropersubgroupG}, $G_{B_n}$ is generated by all $w_0(S(B_k))$ ($1 \leq k \leq n$).
On the other hand, by Lemma \ref{lem:defofpropersubgroupG} again, the product $G_{D_n}\langle s_1 \rangle$ is a semidirect product with $G_{D_n}$ normal, and it is generated by all $w_0(S(D_k))$ ($1 \leq k \leq n$).\\ \indent
We prove the two claims in parallel.
Let $\mathcal{T}=B$ and $L=1$ (for (i)), $\mathcal{T}=D$ and $L=2$ (for (ii)), respectively.
By the above remark, it is enough to show that the group in the left side is generated by all $w_0(S(\mathcal{T}_k))$ ($1 \leq k \leq n$).
We use induction on $n$.
First, note that $w_0(S(\mathcal{T}_n)) \in N_W(W_{S(\mathcal{T}_i)})$ for all $L \leq i \leq n-1$.
Put $W'=W_{S(\mathcal{T}_{n-1})}$.
Then by Corollary \ref{cor:Nofmaximalparabolic} (ii), we have $N_W(W')=W' \rtimes \langle w_0(S(\mathcal{T}_n)) \rangle$.
Thus the claim holds if $n=L+1$; in fact, in this case, $W'=W_{S(\mathcal{T}_L)}$ is generated by all $w_0(S(\mathcal{T}_i))$ ($1 \leq i \leq L$).\\ \indent
If $n>L+1$, then the above equality implies that
\begin{equation*}
\begin{split}
\bigcap_{i=L}^{n-1}N_W(W_{S(\mathcal{T}_i)})&=\left(\bigcap_{i=L}^{n-2}N_W(W_{S(\mathcal{T}_i)})\right) \cap \Bigl(W' \rtimes \langle w_0(S(\mathcal{T}_n)) \rangle\Bigr)\\
&=\left(\bigcap_{i=L}^{n-2}N_{W'}(W_{S(\mathcal{T}_i)})\right) \rtimes \langle w_0(S(\mathcal{T}_n)) \rangle
\end{split}
\end{equation*}
since $w_0(S(\mathcal{T}_n)) \in \bigcap_{i=L}^{n-2}N_W(W_{S(\mathcal{T}_i)})$.
By the induction, the first factor of the semidirect product is generated by all $w_0(S(\mathcal{T}_i))$ ($1 \leq i \leq n-1$).
Thus the claim also holds in this case.
Hence the proof is concluded.
\end{proof}
On the other hand, we have some more properties of the normalizers, which can be deduced without results in \cite{Bri-How} and \cite{How}.
First, we have:
\begin{align}
\label{eq:intersectionofnormalizer}
&\text{If } I,J \subset S, \text{ then } N_W(W_I) \cap N_W(W_J) \subset N_W(W_{I \cap J}).\\
\label{eq:charofnormalizer}
&\text{ For } I \subset S, w \in N_W(W_I) \text{ if and only if } w \cdot \Phi_I=\Phi_I.
\end{align}
((\ref{eq:intersectionofnormalizer}) follows from the well-known fact $W_I \cap W_J=W_{I \cap J}$.
(\ref{eq:charofnormalizer}) follows immediately from (\ref{eq:reflectioninparabolic}).)
Moreover, we have the following:
\begin{lem}
\label{lem:intersectionofnormalizers}
Let $I \subset J \subset S$ such that $J \smallsetminus I \subset I^\perp$.
Then
\begin{equation*}
N_W(W_J) \cap N_W(W_I) \subset N_W(W_{J \smallsetminus I}).
\end{equation*}
\end{lem}
\begin{proof}
Let $w \in N_W(W_J) \cap N_W(W_I)$ and $s \in J \smallsetminus I$.
Then $w \cdot \Phi_J=\Phi_J$ and $w \cdot \Phi_I=\Phi_I$ by (\ref{eq:charofnormalizer}), so that we have $w \cdot \alpha_s \in \Phi_J$ and $w \cdot \alpha_s \not\in \Phi_I$.
Now by the hypothesis and (\ref{eq:supportisconnected}), we have ${\rm supp}(w \cdot \alpha_s) \subset J \smallsetminus I$ and so $w \cdot \alpha_s \in \Phi_{J \smallsetminus I}$.
Hence the claim follows from (\ref{eq:charofnormalizer}).
\end{proof}
%
\section{Main results}\label{sec:mainresults}
\subsection{Direct indecomposability}\label{sec:indecomposability}
In this subsection, we give the main result of this paper that all infinite irreducible Coxeter groups are in fact directly indecomposable, even if it has infinite rank (Theorem \ref{thm:indecomposabilityofW}).
As is mentioned in Introduction, this result was already shown in \cite{Par} for the case of finite rank, in which the finiteness of the ranks is essential and so cannot be removed immediately.\\ \indent
Our proof is based on the following complete description (proved in later sections) of the centralizers of normal subgroups, which are generated by involutions, in irreducible Coxeter groups (possibly of infinite rank):
\begin{thm}
\label{thm:Zofinvolutivenormalsubgroup}
(cf.\ Definition \ref{defn:Coxetergraphs} for notations.)
Let $(W,S)$ be an irreducible Coxeter system of an arbitrary rank, and $H \lhd W$ a normal subgroup generated by involutions.
Then:\\
{\bf (i)} If $H \subset Z(W)$, then $Z_W(H)=W$.\\
{\bf (ii)} If $(W,S)=(W(B_n),S(B_n))$, $2 \leq n \leq \infty$, $\tau \in {\rm Aut}(\varGamma(B_n))$, $H \not\subset Z(W)$ and $H \subset \tau(G_{B_n})$, then $Z_W(H)=\tau(G_{B_n})$.
(cf.\ Lemma \ref{lem:defofpropersubgroupG} for definition of $G_{B_n}$.)\\
{\bf (iii)} If $(W,S)=(W(D_n),S(D_n))$, $3 \leq n \leq \infty$, $\tau \in {\rm Aut}(\varGamma(D_n))$, $H \not\subset Z(W)$ and $H \subset \tau(G_{D_n})$, then $Z_W(H)=\tau(G_{D_n})$.
(cf.\ Lemma \ref{lem:defofpropersubgroupG} for definition of $G_{D_n}$.)\\
{\bf (iv)} Otherwise, $Z_W(H)=Z(W)$.
\end{thm}
This theorem yields the following corollary.
A group $G$ is said to be a {\it central product} of two subgroups $H_1,H_2$ if $G=H_1H_2$ and $H_2 \subset Z_G(H_1)$ (or equivalently $H_1 \subset Z_G(H_2)$).
Note that $H_1 \cap H_2 \subset Z(G)$ in this case.
\begin{cor}
\label{cor:centralproductofW}
Let $(W,S)$ be an irreducible Coxeter system of an arbitrary rank, and suppose that $W$ is a central product of two subgroups $G_1,G_2$ generated by involutions.
Then either $G_1 \subset Z(W)$ or $G_2 \subset Z(W)$.
\end{cor}
\begin{proof}
By definition, we have $G_2 \subset Z_W(G_1)$, $W=G_1Z_W(G_1)$ and $G_1 \lhd W$.
Now if $G_1$ satisfies the condition of cases (ii) or (iii) of Theorem \ref{thm:Zofinvolutivenormalsubgroup}, then $G_1$ and $Z_W(G_1)$ are contained in the same proper subgroup of $W$.
This is impossible, so that we have $G_1 \subset Z(W)$ (case (i)) or $G_2 \subset Z_W(G_1)=Z(W)$ (case (iv)).
\end{proof}
Now our main result follows immediately:
\begin{thm}
\label{thm:indecomposabilityofW}
The only nontrivial direct product decompositions of an irreducible Coxeter group $W$ (of an arbitrary rank) are the ones given in Theorem \ref{thm:Zisdirectfactor}.
In particular, $W$ is directly indecomposable if and only if $W \not\simeq W(\mathcal{T})$ for $\mathcal{T}=B_{2k+1}$, $I_2(4k+2)$ ($k \geq 1$), $E_7$, $H_3$.
\end{thm}
\begin{proof}
Assume that $W=G_1 \times G_2$ for nontrivial subgroups $G_1,G_2 \subset W$.
Then both $G_1$ and $G_2$ are generated by involutions, since $W$ is so.
Thus by Corollary \ref{cor:centralproductofW}, we have either $G_1=Z(W)$ or $G_2=Z(W)$ (since $G_1,G_2 \neq 1$ and $|Z(W)| \leq 2$).
Hence $Z(W) \neq 1$ and so the claim follows from Theorem \ref{thm:Zisdirectfactor}.
\end{proof}
%
\subsection{The Isomorphism Problem}\label{sec:isomorphismproblem}
By using these results, we give some results on the Isomorphism Problem of general Coxeter groups.
Let $(W,S)$ be a Coxeter system with canonical direct product decomposition $W=\prod_{\omega \in \Omega}W_\omega$ into irreducible components $W_\omega$.
Then we put
\begin{equation*}
\Omega_{\rm fin}=\{\omega \in \Omega \mid |W_\omega|<\infty\},\ \Omega_{\rm inf}=\Omega \smallsetminus \Omega_{\rm fin}, W_{\rm fin}=\prod_{\omega \in \Omega_{\rm fin}}W_\omega,\ W_{\rm inf}=\prod_{\omega \in \Omega_{\rm inf}}W_\omega.
\end{equation*}
(Note that $W=W_{\rm fin} \times W_{\rm inf}$.)
Moreover, we write $\Omega_{\mathcal{T}}=\{\omega \in \Omega \mid W_\omega \simeq W(\mathcal{T})\}$ for any type $\mathcal{T}$.
Now our result (proved later) is as follows:
\begin{thm}
\label{thm:isomorphismproblemforW}
(See notations above.)
Let $(W,S)$, $(W',S')$ be two Coxeter systems with the decompositions $W=\prod_{\omega \in \Omega}W_\omega$, $W'=\prod_{\omega' \in \Omega'}W'_{\omega'}$ into irreducible components.
Let $\pi_\omega:W \to W_\omega$, $\pi'_{\omega'}:W' \to W'_{\omega'}$ denote the projections.\\
{\bf (i)} $W \simeq W'$ if and only if the following two conditions are satisfied:\\ \indent
{\bf (I)} There is a bijection $\varphi:\Omega_{\rm inf} \to \Omega'_{\rm inf}$ such that $W_\omega \simeq W'_{\varphi(\omega)}$ for all $\omega \in \Omega_{\rm inf}$.\\ \indent
{\bf (II)} Each of the following subsets of $\Omega$ has the same cardinality as the corresponding subset of $\Omega'$:
\begin{equation*}
\begin{split}
&\Omega_{A_1} \cup \bigl(\bigcup_{k \geq 1}\Omega_{B_{2k+1}}\bigr) \cup \Omega_{E_7} \cup \Omega_{H_3} \cup \bigl(\bigcup_{k \geq 1}\Omega_{I_2(4k+2)}\bigr),\quad \Omega_{B_3} \cup \Omega_{A_3},\\
&\Omega_{B_{2k+1}} \cup \Omega_{D_{2k+1}},\quad \Omega_{I_2(6)} \cup \Omega_{A_2},\quad \Omega_{I_2(4k+2)} \cup \Omega_{I_2(2k+1)}\quad (k \geq 2),\\
&\Omega_{\mathcal{T}}\quad \text{ for } \mathcal{T}=A_n\ (4 \leq n<\infty),\ B_n\ (n<\infty \text{ even}),\ D_n\ (4 \leq n<\infty \text{ even}),\\
&\hspace*{7em} E_6,\ E_7,\ E_8,\ F_4,\ H_3,\ H_4,\ I_2(4k)\ (2 \leq k<\infty).
\end{split}
\end{equation*}
{\bf (ii)} Suppose that $W \simeq W'$, and let $f \in {\rm Isom}(W,W')$.
Then:\\ \indent
{\bf (I)} $f(W_{\rm fin})=W'_{\rm fin}$ (and so the map $g_{\rm fin}$ defined by $g_{\rm fin}=f|_{W_{\rm fin}}$ is an isomorphism $W_{\rm fin} \to W'_{\rm fin}$).\\ \indent
{\bf (II)} There is a bijection $\varphi:\Omega_{\rm inf} \to \Omega'_{\rm inf}$ such that for all $\omega \in \Omega_{\rm inf}$, the map $g_\omega=\pi'_{\varphi(\omega)} \circ f|_{W_\omega}$ is an isomorphism $W_\omega \to W'_{\varphi(\omega)}$.\\ \indent
{\bf (III)} Moreover, there is a map $g_Z \in {\rm Hom}(W_{\rm inf},Z(W'))$ such that
\begin{equation*}
f(w)=
\begin{cases}
g_\omega(w)g_Z(w) & \text{ if } \omega \in \Omega_{\rm inf},\ w \in W_\omega,\\
g_{\rm fin}(w) & \text{ if } w \in W_{\rm fin}.
\end{cases}
\end{equation*}
\end{thm}
Note that this is an analogue of the Krull-Remak-Schmidt Theorem on direct product decompositions of groups, and follows from that (and Theorem \ref{thm:indecomposabilityofW}) if $W$ has a composition series.
(More precisely, the key property in the proof of the K-R-S Theorem, which follows from the existence of composition series, is that any surjective normal endomorphism of an indecomposable factor is either nilpotent or isomorphic.
However, it is not clear whether or not an irreducible Coxeter group has this property.)
Our result here is also a generalization of a result of \cite{Par}.\\ \indent
In order to prove this theorem, we introduce the following ``modified version'' of irreducible components.
Here a group $G$ is said to be {\it admissible} if either $G$ is a nontrivial directly indecomposable irreducible Coxeter group (cf.\ Theorem \ref{thm:indecomposabilityofW}) or $G$ is isomorphic to one of $W(E_7)^+$, $W(H_3)^+$.
\begin{rem}
\label{rem:admissibledecompositionofW}
Let $W=\prod_{\omega \in \Omega}W_\omega$ be the usual decomposition of a Coxeter group $W$ into irreducible components.
Then, by subdividing every directly decomposable $W_\omega$ into the direct factors (cf.\ Theorem \ref{thm:indecomposabilityofW}), we can obtain another decomposition $W=\prod_{\lambda \in \Lambda}G_\lambda$ into admissible subgroups $G_\lambda$.
Moreover, since any infinite $W_\omega$ is directly indecomposable, we can take the index set $\Lambda$ so that $\Omega_{\rm inf} \subset \Lambda$ and $G_\omega=W_\omega$ for all $\omega \in \Omega_{\rm inf}$.
\end{rem}
From now, we consider a family $\mathcal{G}$ of groups which includes all the components of given direct product decompositions.
In our argument below, this family $\mathcal{G}$ is assumed to satisfy the following conditions:
\begin{equation}
\label{eq:conditionforfamily_1}
\begin{split}
&\text{ If } G=\prod_{\lambda \in \Lambda}G_\lambda,\ G_\lambda,G' \in \mathcal{G} \text{ ($\lambda \in \Lambda$) and } f \in {\rm Hom}(G,G') \text{ is surjective},\\
&\text{ then } f \text{ maps a } G_\lambda \text{ onto } G' \text{ (so that it maps all other } G_\mu \text{ into } Z(G')\text{)}.
\end{split}
\end{equation}
\begin{equation}
\label{eq:conditionforfamily_2}
\text{ If } G \in \mathcal{G}, \text{ then } Z(G)=1 \text{ or } Z(G) \text{ is a cyclic group of prime order}.
\end{equation}
(Actually, the condition (\ref{eq:conditionforfamily_2}) can be slightly weakened to the form that $Z(G)$ is either trivial or a finite elementary abelian $p$-group with $p$ prime.
But we omit the detail here, since we do not need such a generalization in this paper.)
\begin{rem}
\label{rem:remarkonthefamilyG}
{\bf (i)} If $\mathcal{G}$ satisfies (\ref{eq:conditionforfamily_1}), then all groups $G \in \mathcal{G}$ are directly indecomposable.
In fact, if $G$ admits a nontrivial decomposition $G=G_1 \times G_2$ with projections $\pi_i:G \to G_i$ ($i=1,2$), then the map $G \times G \to G$, $(w,u) \mapsto \pi_1(w)\pi_2(u)$ is surjective but does not satisfy the conclusion of (\ref{eq:conditionforfamily_1}).\\
{\bf (ii)} If $\mathcal{G}$ satisfies (\ref{eq:conditionforfamily_1}) and (\ref{eq:conditionforfamily_2}), then any $G \in \mathcal{G}$ has the three properties (I)--(III) in Lemma \ref{lem:conditionforvanishingcenter} whenever $Z(G) \neq G$.
This follows immediately from (i).
\end{rem}
\begin{lem}
\label{lem:hombetweenadmissible}
Any family $\mathcal{G}$ of admissible groups satisfies the two conditions.
\end{lem}
\begin{proof}
The condition (\ref{eq:conditionforfamily_2}) follows from Lemma \ref{lem:evensubgroup}.
For (\ref{eq:conditionforfamily_1}), we may assume $G' \not\simeq W(A_1)$ (so that $Z(G') \neq G'$), since otherwise the conclusion is obvious.
Then there is an index $\lambda \in \Lambda$ such that $f(G_\lambda) \not\subset Z(G')$.
Put $G_1=G_\lambda$ and $G_2=\prod_{\mu \in \Lambda \smallsetminus \{\lambda\}}G_\mu$.
Then the hypothesis of (\ref{eq:conditionforfamily_1}) implies that $G'$ is a central product (cf.\ Section \ref{sec:indecomposability}) of $f(G_1)$ and $f(G_2)$, so that $f(G_1) \cap f(G_2) \subset Z(G')$.
Thus the conclusion follows from Lemma \ref{lem:evensubgroup} if $G' \simeq W(E_7)^+$ or $W(H_3)^+$ (in fact, the central product is a direct product since $Z(G')=1$, while $G'$ is directly indecomposable).\\ \indent
On the other hand, suppose that $G'$ is a directly indecomposable irreducible Coxeter group.
Since both $G_1$ and $G_2$ are generated by involutions (cf.\ Lemma \ref{lem:evensubgroup}), $f(G_1)$ and $f(G_2)$ also have this property.
Thus we have $f(G_2) \subset Z(G')$ by Corollary \ref{cor:centralproductofW} (since $f(G_1) \not\subset Z(G')$).
Now if $Z(G') \not\subset f(G_1)$ (so that $f(G_1) \cap Z(G')=1$ since $|Z(G')| \leq 2$), then the central product becomes a (nontrivial) direct product, but this is impossible.
This implies that $f(G_2) \subset Z(G') \subset f(G_1)$ and so $f(G_1)=G'$.
Hence the claim holds.
\end{proof}
\begin{rem}
By a similar argument, it is deduced that any family $\mathcal{G}$, consisting of cyclic groups of prime order and directly indecomposable groups with trivial center, also satisfies the conditions (\ref{eq:conditionforfamily_1}) and (\ref{eq:conditionforfamily_2}).
\end{rem}
We prepare some more notations.
For a decomposition $G=\prod_{\lambda \in \Lambda}G_\lambda$ of $G$, put
\begin{equation}
\label{eq:notationforadmissibledecomposition}
\begin{split}
&G_{\Lambda'}=\prod_{\lambda \in \Lambda'}G_\lambda \text{ (for $\Lambda' \subset \Lambda$)},\ \Lambda_Z=\{\lambda \mid Z(G_\lambda)=G_\lambda\},\ \Lambda_{\neg Z}=\Lambda \smallsetminus \Lambda_Z,\\
&\Lambda_p=\{\lambda \mid |Z(G_\lambda)|=p\}, \Lambda_{Z,p}=\Lambda_Z \cap \Lambda_p, \Lambda_{\neg Z,p}=\Lambda_{\neg Z} \cap \Lambda_p \text{ ($p$ prime or $1$)}.
\end{split}
\end{equation}
Note that the proof of the following theorem is essentially the same as the proof of Theorem 2.1 of \cite{Par}, but slightly more delicate by lack of the assumption on finiteness of the index sets (not only by generality of the context).
Note also that this is also an analogue of the Krull-Remak-Schmidt Theorem.
\begin{thm}
\label{thm:isomorphismproblem}
(See notations above.)
Let $G=\prod_{\lambda \in \Lambda}G_\lambda$, $G'=\prod_{\lambda' \in \Lambda'}G'_{\lambda'}$ be decompositions of two groups $G$, $G'$ into nontrivial subgroups.
Let $\pi_\lambda:G \to G_\lambda$ and $\pi'_{\lambda'}:G' \to G'_{\lambda'}$ be the projections.
Suppose that $\mathcal{G}=\{G_\lambda \mid \lambda \in \Lambda\} \cup \{G'_{\lambda'} \mid \lambda' \in \Lambda'\}$ satisfies the conditions (\ref{eq:conditionforfamily_1}) and (\ref{eq:conditionforfamily_2}).
Let $f \in {\rm Isom}(G,G')$.
Then:\\
{\bf (i)} There is a bijection $\varphi:\Lambda \to \Lambda'$ such that $G_\lambda \simeq G'_{\varphi(\lambda)}$ for all $\lambda \in \Lambda$.
Moreover, for any $\lambda \in \Lambda_{\neg Z}$, the map $g_\lambda=\pi'_{\varphi(\lambda)} \circ f|_{G_\lambda}$ is an isomorphism $G_\lambda \to G'_{\varphi(\lambda)}$.\\
{\bf (ii)} Moreover, there is a map $g_Z \in {\rm Hom}(G,Z(G'))$ such that
\begin{equation*}
f(w)=
\begin{cases}
g_\lambda(w)g_Z(w) & \text{ if } \lambda \in \Lambda_{\neg Z},\ w \in G_\lambda,\\
g_Z(w) & \text{ if } w \in G_{\Lambda_Z}
\end{cases}
\end{equation*}
and that $\pi'_{\varphi(\lambda)} \circ g_Z(G_\lambda)=1$ for all $\lambda \in \Lambda_{\neg Z}$.\\
{\bf (iii)} If $\bigcup_{p \neq 1}\Lambda_p \subset \Lambda^\natural \subset \Lambda$, then $\bigcup_{p \neq 1}\Lambda'_p \subset \varphi(\Lambda^\natural)$ and $f(G_{\Lambda^\natural})=G'_{\varphi(\Lambda^\natural)}$.
\end{thm}
\begin{proof}
Note that $\bigcup_{p \neq 1}\Lambda_p=\{\lambda \in \Lambda \mid Z(G_\lambda) \neq 1\}$.
Then the claim (iii) is deduced from the other claims (since now $Z(G) \subset G_{\Lambda^\natural}$ and $Z(G') \subset G'_{\varphi(\Lambda^\natural)}$).\\ \indent
From now, we prove the claims (i) and (ii).
First, we put (symmetrically)
\begin{equation*}
f_{\lambda'}=\pi'_{\lambda'} \circ f \in {\rm Hom}(G,G'_{\lambda'}) \text{ ($\lambda' \in \Lambda'$)},\ f'_\lambda=\pi_\lambda \circ f^{-1} \in {\rm Hom}(G',G_\lambda) \text{ ($\lambda \in \Lambda$)},
\end{equation*}
and define (symmetrically)
\begin{equation*}
\begin{split}
\mathcal{A}'_\lambda=\{\lambda' \in \Lambda' \mid f_{\lambda'}(G_\lambda) \not\subset Z(G'_{\lambda'})\} \subset \Lambda'_{\neg Z} \text{ for } \lambda \in \Lambda_{\neg Z},\\
\mathcal{A}_{\lambda'}=\{\lambda \in \Lambda \mid f'_\lambda(G'_{\lambda'}) \not\subset Z(G_\lambda)\} \subset \Lambda_{\neg Z} \text{ for } \lambda' \in \Lambda'_{\neg Z}.
\end{split}
\end{equation*}
Note that $\mathcal{A}'_\lambda \neq \emptyset$ since $f(G_\lambda) \not\subset Z(G')$ (and $\mathcal{A}_{\lambda'} \neq \emptyset$ by symmetry).
Moreover, since $f_{\lambda'}:G \to G'_{\lambda'}$ is surjective, the condition (\ref{eq:conditionforfamily_1}) implies that
\begin{equation*}
\text{ if } \lambda' \in \mathcal{A}'_\lambda, \text{ then } f_{\lambda'}(G_\lambda)=G'_{\lambda'} \text{ and } f_{\lambda'}(G_\mu) \subset Z(G'_{\lambda'}) \text{ for all } \mu \in \Lambda \smallsetminus \{\lambda\}.
\end{equation*}
By symmetry, a similar property holds for $\lambda \in \mathcal{A}_{\lambda'}$ (with respect to the map $f'_\lambda$).\\ \indent
We prove the following claims:\\ \indent
{\bf Claim 1:}\ If $\lambda,\mu \in \Lambda_{\neg Z}$ and $\lambda \neq \mu$, then $\mathcal{A}'_\lambda \cap \mathcal{A}'_\mu=\emptyset$.\\ \indent
{\bf Claim 2:}\ If $\lambda' \in \mathcal{A}'_\lambda$, then $\lambda \in \mathcal{A}_{\lambda'}$.
(Thus $|\mathcal{A}'_\lambda|=1$ for all $\lambda \in \Lambda_{\neg Z}$, by Claim 1 and symmetry.
Moreover, by symmetry, the map $\varphi:\Lambda_{\neg Z} \to \Lambda'_{\neg Z}$ defined by $\mathcal{A}'_\lambda=\{\varphi(\lambda)\}$ is a bijection with inverse map satisfying $\mathcal{A}_{\lambda'}=\{\varphi^{-1}(\lambda')\}$.)\\ \indent
{\bf Claim 3:}\ The map $g_\lambda$ ($\lambda \in \Lambda_{\neg Z}$) in (i) is an isomorphism $G_\lambda \to G'_{\varphi(\lambda)}$.\\ \indent
{\bf Claim 4:}\ $f(Z(G_{\Lambda_{\neg Z,p}}))=Z(G'_{\Lambda'_{\neg Z,p}})$ for all primes $p$.\\ \indent
{\bf Claim 5:}\ For each prime $p$, $\Lambda_{Z,p}$ and $\Lambda'_{Z,p}$ have the same cardinality.\\ \indent
{\bf Proof of Claim 1:}\ 
Assume contrary that $\lambda' \in \mathcal{A}'_\lambda \cap \mathcal{A}'_\mu$.
Then the relation $\lambda' \in \mathcal{A}'_\lambda$ means that $f_{\lambda'}(G_\lambda) \not\subset Z(G'_{\lambda'})$, while the relation $\lambda' \in \mathcal{A}'_\mu$ implies (by the above property) that $f_{\lambda'}(G_\lambda) \subset Z(G'_{\lambda'})$ (since $\lambda \neq \mu$).
This is a contradiction.\\ \indent
{\bf Proof of Claim 2:}\ 
Since $G'_{\lambda'} \neq Z(G'_{\lambda'})$, we can take an element $w \in G'_{\lambda'} \smallsetminus Z(G'_{\lambda'})$.
Put $u_\mu=f'_\mu(w) \in G_\mu$ for $\mu \in \Lambda$, so that we have $w=f(\prod_{\mu \in \Lambda}u_\mu)$.
Now $f_{\lambda'}(u_\mu) \in Z(G'_{\lambda'})$ for all $\mu \in \Lambda \smallsetminus \{\lambda\}$, while $w=\pi'_{\lambda'}(w) \not\in Z(G'_{\lambda'})$.
Thus we have $f_{\lambda'}(u_\lambda) \not\in Z(G'_{\lambda'})$ and so $u_\lambda \not\in Z(G_\lambda)$ (since $f_{\lambda'}(G_\lambda)=G'_{\lambda'}$).
Hence $\lambda \in \mathcal{A}_{\lambda'}$.\\ \indent
{\bf Proof of Claim 3:}\ 
Note that $g_\lambda:G_\lambda \to G'_{\varphi(\lambda)}$ is surjective (as above).
Now the following equivalence holds for all $w \in G_\lambda$:
\begin{equation*}
f_{\varphi(\lambda)}(w) \in Z(G'_{\varphi(\lambda)}) \Longleftrightarrow f(w) \in Z(G') \Longleftrightarrow w \in Z(G) \Longleftrightarrow w \in Z(G_\lambda)
\end{equation*}
(we use the fact $\mathcal{A}'_\lambda=\{\varphi(\lambda)\}$ for the first equivalence).
This implies that $\ker g_\lambda$ is contained in the simple group $Z(G_\lambda)$ (cf.\ (\ref{eq:conditionforfamily_2})), so that $\ker g_\lambda=1$ or $Z(G_\lambda)$.
Thus $g_\lambda$ is injective (and so an isomorphism) if $Z(G_\lambda)=1$.
Moreover, if $Z(G'_{\varphi(\lambda)})=1$, then $f'_\lambda|_{G'_{\varphi(\lambda)}}$ is an isomorphism $G'_{\varphi(\lambda)} \to G_\lambda$ by symmetry, so that we have $Z(G_\lambda)=1$.
Thus $g_\lambda$ is injective (as above) also in this case.\\ \indent
On the other hand, suppose $Z(G'_{\varphi(\lambda)}) \neq 1$.
Then by the above equivalence, there is an element $w \in Z(G_\lambda)$ such that $g_\lambda(w) \neq 1$ (since $g_\lambda$ is surjective).
Thus we have $\ker g_\lambda \neq Z(G_\lambda)$ and so $\ker g_\lambda=1$.
Hence $g_\lambda$ is an isomorphism.\\ \indent
{\bf Proof of Claim 4:}\ 
Note that $Z(G)=\prod_{p \neq 1}Z(G_{\Lambda_p})$ and each $Z(G_{\Lambda_p})$ is an elementary abelian $p$-group, by (\ref{eq:conditionforfamily_2}).
$Z(G')$ also admits a similar decomposition.
Thus the isomorphism $f|_{Z(G)}:Z(G) \to Z(G')$ maps each $Z(G_{\Lambda_p})$ onto $Z(G'_{\Lambda'_p})$.
Moreover, for any $\lambda \in \Lambda_{\neg Z,p}$, the composite homomorphism $G_\lambda \overset{f}{\to} G' \to G'_{\Lambda'_{Z,p}}$ (where the latter map is the projection) maps $Z(G_\lambda)$ to $1$, by Remark \ref{rem:remarkonthefamilyG} (ii) (note that $Z(G'_{\Lambda'_{Z,p}})=G'_{\Lambda'_{Z,p}}$).
Thus we have $f(Z(G_\lambda)) \subset G'_{\Lambda'_{\neg Z,p}}$ for any $\lambda \in \Lambda_{\neg Z,p}$ and so $f(Z(G_{\Lambda_{\neg Z,p}})) \subset Z(G'_{\Lambda'_{\neg Z,p}})$.
Now this claim holds by symmetry.\\ \indent
{\bf Proof of Claim 5:}\ 
Note that $Z(G_{\Lambda_p})=G_{\Lambda_{Z,p}} \times Z(G_{\Lambda_{\neg Z,p}})$ and $Z(G'_{\Lambda'_p})$ admits a similar decomposition.
Moreover, we have $f(Z(G_{\Lambda_p}))=Z(G'_{\Lambda'_p})$ and $f(Z(G_{\Lambda_{\neg Z,p}}))=Z(G'_{\Lambda'_{\neg Z,p}})$ by Claim 4.
Thus the complementary factors $G_{\Lambda_{Z,p}}$, $G'_{\Lambda'_{Z,p}}$, which
are elementary abelian $p$-groups with basis having the same cardinality as $\Lambda_{Z,p}$, $\Lambda'_{Z,p}$ respectively, are also isomorphic.
Now this claim follows from uniqueness of the dimension of a vector space.\\ \indent
{\bf Conclusion.}\ 
Since $\Lambda_Z$, $\Lambda'_Z$ are disjoint unions of $\Lambda_{Z,p}$, $\Lambda'_{Z,p}$ respectively (cf.\ (\ref{eq:conditionforfamily_2})), Claim 5 implies that this $\varphi$ extends (not uniquely) to a bijection $\varphi:\Lambda \to \Lambda'$ satisfying (i) (note that $\Lambda_{Z,1}=\Lambda'_{Z,1}=\emptyset$ by the hypothesis).
Moreover, define a map $g_Z:G \to Z(G')$ componentwise by 
\begin{equation*}
g_Z(w)=
\begin{cases}
\prod_{\lambda' \in \Lambda' \smallsetminus \{\varphi(\lambda)\}}f_{\lambda'}(w) & \text{ if } \lambda \in \Lambda_{\neg Z}, w \in G_\lambda,\\
f(w) & \text{ if } w \in G_{\Lambda_Z}.
\end{cases}
\end{equation*}
Note that $G_{\Lambda_Z} \subset Z(G)$, while in the above definition, we have $f_{\lambda'}(w) \in Z(G'_{\lambda'})$ by the fact $\mathcal{A}'_\lambda=\{\varphi(\lambda)\}$.
Since $Z(G')$ is abelian, these facts imply that $g_Z$ is a well-defined group homomorphism.
Now the claim (ii) follows from definition.
\end{proof}
\begin{proof}[Proof of Theorem \ref{thm:isomorphismproblemforW}]
Let $W=\prod_{\lambda \in \Lambda}G_\lambda$, $W'=\prod_{\lambda' \in \Lambda'}G'_{\lambda'}$ be the decompositions into admissible groups given in Remark \ref{rem:admissibledecompositionofW}.\\
{\bf (i)} Each of the sets in the condition (II), except $\Omega_{E_7}$ and $\Omega_{H_3}$ in the last row, has the same cardinality as the set $\{\lambda \in \Lambda \mid G_\lambda \simeq W(\mathcal{T}')\}$ where $\mathcal{T}'=A_1$, $A_3$, $D_{2k+1}$, $A_2$, $I_2(2k+1)$ and $\mathcal{T}$, respectively (note that no two admissible finite groups of distinct types are isomorphic; cf.\ Lemma \ref{lem:evensubgroup}).
Moreover, each of $\Omega_{E_7}$ and $\Omega_{H_3}$ has the same cardinality as $\{\lambda \in \Lambda \mid G_\lambda \simeq W(\mathcal{T}')^+\}$ for $\mathcal{T}'=E_7$ and $H_3$, respectively.
Similar relations also hold for $W'$.
Thus the two conditions (I), (II) are satisfied if and only if there is a bijection $\psi:\Lambda \to \Lambda'$ such that $G_\lambda \simeq G'_{\psi(\lambda)}$ for all $\lambda \in \Lambda$.
Hence the claim follows from Theorem \ref{thm:isomorphismproblem} (i) (which can be applied indeed to the case, by Lemma \ref{lem:hombetweenadmissible}).\\
{\bf (ii)} Take $\varphi:\Lambda \to \Lambda'$, $g_\lambda \in {\rm Isom}(G_\lambda,G'_{\varphi(\lambda)})$ ($\lambda \in \Lambda_{\neg Z}$) and $g'_Z \in {\rm Hom}(W,Z(W'))$ as in the conclusion of Theorem \ref{thm:isomorphismproblem}.
By Remark \ref{rem:admissibledecompositionofW}, $g_\omega \in {\rm Isom}(W_\omega,W'_{\varphi(\omega)})$ for all $\omega \in \Omega_{\rm inf}$, so that the claim (II) holds.
The claim (I) follows from Theorem \ref{thm:isomorphismproblem} (iii) (by putting $\Lambda^\natural=\Lambda \smallsetminus \Omega_{\rm inf}$).
Moreover, the claim (III) also follows from Theorem \ref{thm:isomorphismproblem}, by putting $g_Z=g'_Z|_{W_{\rm inf}}$.
Hence the proof is concluded.
\end{proof}
%
\subsection{Automorphism groups}\label{sec:automorphismgroup}
Owing to Theorems \ref{thm:isomorphismproblemforW} and \ref{thm:isomorphismproblem}, we can examine the automorphism groups of $W=\prod_{\omega \in \Omega}W_\omega$ and $G=\prod_{\lambda \in \Lambda}G_\lambda$ respectively (Theorem \ref{thm:automorphismgroup}), under the hypothesis in Section \ref{sec:isomorphismproblem}.
In this subsection, the {\it complete} direct product of groups is denoted by a symbol $\overline{\prod}$.\\ \indent
As is remarked in Section \ref{sec:generalgroups}, if $G',G''$ are groups and $G'=G'_1 \times G'_2$, then the set ${\rm Hom}(G'_1,G'')$ is embedded naturally into ${\rm Hom}(G',G'')$.
In this manner, each ${\rm Aut}(G_\lambda)$, ${\rm Aut}(W_\omega)$ is embedded into ${\rm Aut}(G)$, ${\rm Aut}(W)$ respectively.
The group ${\rm Aut}(W_{\rm fin})$ is also embedded into ${\rm Aut}(W)$.\\ \indent
On the other hand, the symmetric group on each isomorphism class of components of $G$ or $W$ is also embedded into the automorphism group, as follows.
For the case of $G$, we partition the index set $\Lambda_{\neg Z}$ into subsets $\Lambda_\xi$ ($\xi \in \Xi$) so that $\lambda,\lambda' \in \Lambda_{\neg Z}$ are in the same subset if and only if $G_\lambda \simeq G_{\lambda'}$.
Moreover, for $\xi \in \Xi$, we choose an ``identity map'' ${\rm id}_{\mu,\lambda} \in {\rm Isom}(G_\lambda,G_\mu)$ for each $\lambda,\mu \in \Lambda_\xi$ so that ${\rm id}_{\lambda,\lambda}={\rm id}_{G_\lambda}$, ${\rm id}_{\lambda,\mu}={{\rm id}_{\mu,\lambda}}^{-1}$ and ${\rm id}_{\nu,\mu} \circ {\rm id}_{\mu,\lambda}={\rm id}_{\nu,\lambda}$ for all $\lambda,\mu,\nu \in \Lambda_\xi$.
(This can be done by taking a maximal tree in the category of groups $G_\lambda$ ($\lambda \in \Lambda_\xi$) and group isomorphisms.)
Then each element $\tau$ of the symmetric group ${\rm Sym}(\Lambda_\xi)$ on $\Lambda_\xi$ induces an automorphism of the factor $G_{\Lambda_\xi}$ of $G$; namely, 
\begin{equation*}
\tau(w)={\rm id}_{\tau(\lambda),\lambda}(w) \in G_{\tau(\lambda)} \text{ for } \lambda \in \Lambda_\xi \text{ and } w \in G_\lambda.
\end{equation*}
In this manner, ${\rm Sym}(\Lambda_\xi)$ is embedded into ${\rm Aut}(G_{\Lambda_\xi})$, and so also into ${\rm Aut}(G)$.
Similarly, we write $\Omega=\bigsqcup_{\upsilon \in \Upsilon}\Omega_\upsilon$, choose ``identity maps'' ${\rm id}_{\omega',\omega} \in {\rm Isom}(W_\omega,W_{\omega'})$ and then embed every symmetric group ${\rm Sym}(\Omega_\upsilon)$ into ${\rm Aut}(W)$.
Moreover, put
\begin{equation*}
\Upsilon_{\rm fin}=\{\upsilon \in \Upsilon \mid |W_\omega|<\infty \text{ for } \omega \in \Omega_\upsilon\} \text{ and } \Upsilon_{\rm inf}=\Upsilon \smallsetminus \Upsilon_{\rm fin}.
\end{equation*}
\indent
For a group $G'$, recall (Section \ref{sec:generalgroups}) the structure of the monoid ${\rm Hom}(G',Z(G'))$, the action of ${\rm Aut}(G')$ on it and the embedding $f \mapsto f^\flat$ into the monoid ${\rm End}(G')$ compatible with the action of ${\rm Aut}(G')$.
By this map, the group ${\rm Hom}(G',Z(G'))^\times$ of invertible elements of ${\rm Hom}(G',Z(G'))$ is embedded into ${\rm Aut}(G')$.\\ \indent
Now for the group $G$, let
\begin{equation*}
\begin{split}
 {\rm Hom}(G,Z(G))_o=\{f \in {\rm Hom}(G,Z(G)) \mid &f(G_{\Lambda_Z})=1,\\
&f(G_\lambda) \subset Z(G_\lambda) \text{ for all } \lambda \in \Lambda_{\neg Z}\}
\end{split}
\end{equation*}
(cf.\ (\ref{eq:notationforadmissibledecomposition}) for notations).
Since we assumed that each $G_\lambda$ ($\lambda \in \Lambda_{\neg Z}$) satisfies the three conditions in Lemma \ref{lem:conditionforvanishingcenter} (cf.\ Remark \ref{rem:remarkonthefamilyG} (ii)), we have $f(Z(G))=1$ for all $f \in {\rm Hom}(G,Z(G))_o$.
Thus by Lemma \ref{lem:subgroupofHom} (i), ${\rm Hom}(G,Z(G))_o$ is an abelian subgroup of ${\rm Hom}(G,Z(G))^\times$ with multiplication $(f*g)(w)=f(w)g(w)$ ($f,g \in {\rm Hom}(G,Z(G))_o$, $w \in G$).\\ \indent
On the other hand, since $Z(W_{\rm inf})=1$, Lemma \ref{lem:subgroupofHom} (ii) implies that the set ${\rm Hom}(W_{\rm inf},Z(W))$ forms an abelian normal subgroup of ${\rm Hom}(W,Z(W))^\times$ with multiplication $(f*g)(w)=f(w)g(w)$ ($f,g \in {\rm Hom}(W_{\rm inf},Z(W))$, $w \in W_{\rm inf}$).
Since now $Z(W)$ is an elementary abelian $2$-group, ${\rm Hom}(W_{\rm inf},Z(W))$ is also an elementary abelian $2$-group.\\ \indent
Now our result is stated as follows:
\begin{thm}
\label{thm:automorphismgroup}
(See notations above.)\\
{\bf (i)} Put $H_1={{\rm Hom}(G,Z(G))^\times}^\flat$, $H_2=\overline{\prod}_{\lambda \in \Lambda_{\neg Z}}{\rm Aut}(G_\lambda)$, $H_3=\overline{\prod}_{\xi \in \Xi}{\rm Sym}(\Lambda_\xi)$ and $H_4={\rm Hom}(G,Z(G))_o^\flat$.
Then 
\begin{equation*}
{\rm Aut}(G)=(H_1H_2) \rtimes H_3,\ H_1 \lhd {\rm Aut}(G),\ H_2 \lhd H_2H_3,\ H_1 \cap H_2=H_4.
\end{equation*}
 {\bf (ii)} Put $H'_1={\rm Hom}(W_{\rm inf},Z(W))^\flat$, $H'_2={\rm Aut}(W_{\rm fin})$, $H'_3=\overline{\prod}_{\omega \in \Omega_{\rm inf}}{\rm Aut}(W_\omega)$ and $H'_4=\overline{\prod}_{\upsilon \in \Upsilon_{\rm inf}}{\rm Sym}(\Omega_\upsilon)$.
Then
\begin{equation*}
{\rm Aut}(W)=H'_1 \rtimes (H'_2 \times H'_3) \rtimes H'_4,\ H'_2H'_4=H'_2 \times H'_4,\ H'_3H'_4=H'_3 \rtimes H'_4.
\end{equation*}
{\bf (iii)} The subgroup $H=\left(\overline{\prod}_{\omega \in \Omega}{\rm Aut}(W_\omega)\right)\left(\overline{\prod}_{\upsilon \in \Upsilon}{\rm Sym}(\Omega_\upsilon)\right)$ has finite index in ${\rm Aut}(W)$ if and only if, either $Z(W)=1$ or the odd Coxeter graph (cf.\ Definition \ref{defn:oddCoxetergraph}) $\varGamma^{\rm odd}$ of $W$ consists of only finitely many connected components.
(Hence the index is finite whenever $W$ has finite rank.)
\end{thm}
From now, we prove this theorem.
First, we prove (i) and (ii).
Note that $H'_2H'_3=H'_2 \times H'_3$ and $H'_2H'_4=H'_2 \times H'_4$ by definition.
Moreover, by definition,
\begin{equation}
\label{eq:charofproductofaut}
\begin{split}
&H_2=\{f \in {\rm Aut}(G) \mid f(w)=w \text{ ($w \in G_{\Lambda_Z}$)},\ f(G_\lambda)=G_\lambda \text{ ($\lambda \in \Lambda_{\neg Z}$)}\},\\
&H'_3=\{f \in {\rm Aut}(W) \mid f(w)=w \text{ ($w \in W_{\rm fin}$)},\ f(W_\omega)=W_\omega \text{ ($\omega \in \Omega_{\rm inf}$)}\}.
\end{split}
\end{equation}
{\bf Claim 1.}
{\bf (i)} ${\rm Aut}(G)=H_1H_2H_3$.
{\bf (ii)} ${\rm Aut}(W)=H'_1H'_2H'_3H'_4$.
\begin{proof}
{\bf (i)} Let $f \in {\rm Aut}(G)$, and take $\varphi$, $g_\lambda$, $g_Z$ as in Theorem \ref{thm:isomorphismproblem}.
Note that $\varphi(\Lambda_\xi)=\Lambda_\xi$ for all $\xi \in \Xi$.
Now define $f_1 \in {\rm Hom}(G,Z(G))$ by
\begin{equation*}
f_1(w)=
\begin{cases}
g_Z \circ g_{\varphi^{-1}(\lambda)}^{-1}(w)^{-1} & \text{ for } \lambda \in \Lambda_{\neg Z},\ w \in G_\lambda,\\
wf(w)^{-1} & \text{ for } w \in G_{\Lambda_Z}
\end{cases}
\end{equation*} 
(this is well defined since $G_{\Lambda_Z} \subset Z(G)$).
Then by definition and Theorem \ref{thm:isomorphismproblem}, we have $f={f_1}^\flat \circ f_2 \circ f_3$, where
\begin{equation*}
f_2=(g_{\varphi^{-1}(\lambda)} \circ {\rm id}_{\varphi^{-1}(\lambda),\lambda})_{\lambda \in \Lambda_{\neg Z}} \in H_2,\ f_3=(\varphi|_{\Lambda_\xi})_{\xi \in \Xi} \in H_3.
\end{equation*}
Moreover, we have ${f_1}^\flat=f \circ f_3^{-1} \circ f_2^{-1} \in {\rm Aut}(G)$ and so $f_1 \in {\rm Hom}(G,Z(G))^\times$ by Lemma \ref{lem:HomtoEndisembedding} (ii).
Hence ${f_1}^\flat \in H_1$ and so $f \in H_1H_2H_3$.\\
{\bf (ii)} Let $f \in {\rm Aut}(W)$, and take $\varphi$, $g_{\rm fin}$, $g_\lambda$, $g_Z$ as in Theorem \ref{thm:isomorphismproblemforW} (ii).
Note that $\varphi(\Omega_\upsilon)=\Omega_\upsilon$ for all $\upsilon \in \Upsilon$.
Now define $f_1 \in {\rm Hom}(W_{\rm inf},Z(W))$ by
\begin{equation*}
f_1(w)=g_Z \circ g_{\varphi^{-1}(\omega)}^{-1}(w)^{-1} \text{ for } \omega \in \Omega_{\rm inf},\ w \in W_\omega.
\end{equation*} 
Then we have (by definition and Theorem \ref{thm:isomorphismproblemforW} (ii))
\begin{equation*}
f={f_1}^\flat \circ g_{\rm fin} \circ (g_{\varphi^{-1}(\omega)} \circ {\rm id}_{\varphi^{-1}(\omega),\omega})_{\omega \in \Omega_{\rm inf}} \circ (\varphi|_{\Omega_\upsilon})_{\upsilon \in \Upsilon_{\rm inf}} \in H'_1H'_2H'_3H'_4.
\end{equation*}
Hence the proof is concluded.
\end{proof}
\noindent
{\bf Claim 2.}
{\bf (i)} If $f^\flat \in H_1$, $\lambda,\mu \in \Lambda_{\neg Z}$ and $f^\flat(G_\lambda) \subset G_\mu$, then $\lambda=\mu$ and $f(G_\lambda) \subset Z(G_\lambda)$.\\
{\bf (ii)} If $f^\flat \in H'_1$, $\omega,\omega' \in \Omega_{\rm inf}$ and $f^\flat(W_\omega) \subset W_{\omega'}$, then $\omega=\omega'$ and $f(W_\omega)=1$.
\begin{proof}
{\bf (i)} By the choice of $\lambda$, we can take $w \in G_\lambda \smallsetminus Z(G_\lambda)$.
Now we have $\pi_\lambda(f(w)) \in Z(G_\lambda)$ (where $\pi_\lambda$ is the projection $G \to G_\lambda$) and so $\pi_\lambda(f^\flat(w))=w\pi_\lambda(f(w))^{-1} \neq 1$.
Since $f^\flat(w) \in G_\mu$, this implies that $\mu=\lambda$.
Now the latter part follows from definition of the map $f^\flat$.\\
{\bf (ii)} By a similar argument to (i), we have $\omega=\omega'$ and $f(W_\omega) \subset Z(W_\omega)$.
Hence the claim holds since $Z(W_\omega)=1$.
\end{proof}
\noindent
{\bf Claim 3.}
{\bf (i)} $(H_1H_2) \cap H_3=1$.
{\bf (ii)} $(H'_1H'_2H'_3) \cap H'_4=1$.
\begin{proof}
{\bf (i)} Let $f_1 \in H_1$, $f_2 \in H_2$ such that $f_1 \circ f_2 \in H_3$.
By (\ref{eq:charofproductofaut}) and definition of $H_3$, both $f_2^{-1}$ and $f_1 \circ f_2$ map each component $G_\lambda$ ($\lambda \in \Lambda_{\neg Z}$) onto a component, so that $f_1$ also does so.
By Claim 2 (i), $f_1$ maps each $G_\lambda$ ($\lambda \in \Lambda_{\neg Z}$) onto itself, while $f_2$ also does so (cf.\ (\ref{eq:charofproductofaut})).
Thus $f_1 \circ f_2 \in H_3$ also has this property.
By definition of $H_3$, this occurs only if $f_1 \circ f_2={\rm id}_G$.
Hence the claim holds.\\
{\bf (ii)} The proof is similar to (i); if $f_i \in H'_i$ ($i=1,2,3$) and $f_4=f_1 \circ f_2 \circ f_3 \in H'_4$, then $f_1=f_4 \circ f_3^{-1} \circ f_2^{-1}$ must map each $W_\omega$ ($\omega \in \Omega_{\rm inf}$) onto some component, which is $W_\omega$ by Claim 2 (ii).
This implies that $f_4$ maps each $W_\omega$ ($\omega \in \Omega_{\rm inf}$) onto itself, so that $f_4={\rm id}_W$ by definition of $H'_4$.
Hence the claim holds.
\end{proof}
\noindent
{\bf Claim 4.}
{\bf (i)} $H_2 \lhd H_2H_3$.
{\bf (ii)} $H'_3 \lhd H'_3H'_4$.
\begin{proof}
For (i), it is enough to show that $f_3 \circ f_2 \circ f_3^{-1} \in H_2$ for all $f_2 \in H_2$ and $f_3 \in H_3$.
By definition, $f_3$ is identity on $G_{\Lambda_Z}$ and maps each $G_\lambda$ ($\lambda \in \Lambda_{\neg Z}$) onto a component.
Now by (\ref{eq:charofproductofaut}), $f_3 \circ f_2 \circ f_3^{-1}$ also satisfies the condition in (\ref{eq:charofproductofaut}), so that it belongs to $H_2$.
Hence the claim holds.
The proof of (ii) is similar.
\end{proof}
\noindent
{\bf Claim 5.}
{\bf (i)} $H_1 \lhd {\rm Aut}(G)$.
{\bf (ii)} $H'_1 \lhd {\rm Aut}(W)$.
\begin{proof}
{\bf (i)} Note that ${\rm Aut}(G)$ acts on the monoid ${\rm Hom}(G,Z(G))$.
Thus its subgroup ${\rm Hom}(G,Z(G))^\times$ of the invertible elements is invariant under the action.
Now the claim follows from Lemma \ref{lem:HomtoEndisembedding} (iii).\\
{\bf (ii)} By Lemma \ref{lem:HomtoEndisembedding} (iii), it is enough to show that the subgroup ${\rm Hom}(W_{\rm inf},Z(W))$ of ${\rm Hom}(W,Z(W))$ is invariant under the action of ${\rm Aut}(W)$.
Moreover, by Claim 1, it is enough to show that $h \circ f \circ h^{-1} \in {\rm Hom}(W_{\rm inf},Z(W))$ for all $f \in {\rm Hom}(W_{\rm inf},Z(W))$ and $h \in H'_2H'_3H'_4$.
Now we have $h(W_{\rm fin})=W_{\rm fin}$ by definition of $H'_2$, $H'_3$ and $H'_4$, so that $h \circ f \circ h^{-1}(W_{\rm fin})=h(f(W_{\rm fin}))=h(1)=1$.
Hence the claim holds.
\end{proof}
\noindent
{\bf Claim 6.}
{\bf (i)} $H_1 \cap H_2=H_4$.
{\bf (ii)} $H'_1 \cap (H'_2H'_3)=1$.
\begin{proof}
{\bf (i)} Let $f^\flat \in H_1 \cap H_2$.
Then by (\ref{eq:charofproductofaut}), we have $f^\flat(w)=w$ (or equivalently $f(w)=1$) for all $w \in G_{\Lambda_Z}$ and $f^\flat(G_\lambda)=G_\lambda$ for all $\lambda \in \Lambda_{\neg Z}$.
Thus we have $f \in {\rm Hom}(G,Z(G))_o$ by Claim 2 (i), so that $f^\flat \in H_4$.
Conversely, $H_4 \subset H_1$ by definition, while $H_4 \subset H_2$ by (\ref{eq:charofproductofaut}) and definition of $H_4$.
Hence the claim holds.\\
{\bf (ii)} Let $f^\flat \in H'_1 \cap (H'_2H'_3)$.
Then for any $\omega \in \Omega_{\rm inf}$, we have $f^\flat(W_\omega)=W_\omega$ by definition of $H'_2$ and $H'_3$.
Thus we have $f(W_\omega)=1$ by Claim 2 (ii).
Hence $f=1$ and $f^\flat={\rm id}_W$.
\end{proof}
Now the claims (i) and (ii) of Theorem \ref{thm:automorphismgroup} hold.
Namely:\\
{\bf (i)} We have $H_1 \cap H_2=H_4$ (Claim 6), $H_1 \lhd {\rm Aut}(G)$ (Claim 5), $H_2 \lhd H_2H_3$ (Claim 4) and so ${\rm Aut}(G)=(H_1H_2)H_3$ (Claim 1) $=(H_1H_2) \rtimes H_3$ (Claim 3).\\
{\bf (ii)} We have $H'_2H'_3=H'_2 \times H'_3$, $H'_2H'_4=H'_2 \times H'_4$ (as the above remark), $H'_3H'_4=H'_3 \rtimes H'_4$ (Claims 3, 4) and so ${\rm Aut}(W)=H'_1(H'_2 \times H'_3)H'_4$ (Claim 1) $=\bigl(H'_1(H'_2 \times H'_3)\bigr) \rtimes H'_4$ (Claims 3, 5) $=H'_1 \rtimes (H'_2 \times H'_3) \rtimes H'_4$ (Claims 5, 6).
\begin{proof}[Proof of Theorem \ref{thm:automorphismgroup} (iii)]
If $Z(W)=1$, then all irreducible components of $W$ are directly indecomposable (cf.\ Theorem \ref{thm:indecomposabilityofW}), so that the decomposition $W=\prod_{\omega \in \Omega}W_\omega$ itself satisfies the conditions (\ref{eq:conditionforfamily_1}) and (\ref{eq:conditionforfamily_2}) in Section \ref{sec:isomorphismproblem}.
Thus we can apply the result (i) to this decomposition.
Now $H_1=1$ since $Z(W)=1$.
Moreover, $\Omega=\Omega_{\neg Z}$ in this case, so that we have $H=H_2H_3={\rm Aut}(W)$.\\ \indent
From now, we assume that $Z(W) \neq 1$.
For $f \in {\rm Aut}(W)$, let ${\rm sep}(f)$ be the set of all $\omega \in \Omega$ such that $f(W_\omega) \not\subset W_{\omega'}$ for all $\omega' \in \Omega$.
Since any element of $H$ maps each component $W_\omega$ onto a component, the cardinality of the set ${\rm sep}(f)$ is invariant in each coset of ${\rm Aut}(W)/H$.
Moreover, by definition, we have
\begin{equation*}
H=\left(\overline{\prod}_{\omega \in \Omega_{\rm fin}}{\rm Aut}(W_\omega)\right)\left(\overline{\prod}_{\upsilon \in \Upsilon_{\rm fin}}{\rm Sym}(\Omega_\upsilon)\right) \times H'_3H'_4 \subset H'_2 \times (H'_3H'_4).
\end{equation*}
\indent
{\bf Case 1. $\varGamma^{\rm odd}$ consists of only finitely many connected components:}\ 
This implies that $|\Omega|<\infty$ and $|{\rm Hom}(W_{\rm inf},\{\pm 1\})|<\infty$ (cf.\ Lemma \ref{lem:charofWtopm1}).
Since $Z(W)$ is now a finite elementary abelian $2$-group, (ii) implies that $H'_2H'_3H'_4$ has index $|H'_1|=|{\rm Hom}(W_{\rm inf},Z(W))|<\infty$ in ${\rm Aut}(W)$.
Moreover, since now $|W_{\rm fin}|<\infty$, the index of $H$ in $H'_2H'_3H'_4$ is $\leq |H'_2|<\infty$.
Thus $H$ has finite index also in ${\rm Aut}(W)$.\\ \indent
{\bf Case 2. $\varGamma^{\rm odd}$ consists of infinitely many connected components:}\ 
Now we have to show that $H$ has infinite index in ${\rm Aut}(W)$.\\ \indent
{\bf Subcase 2-1. The odd Coxeter graph of some $W_\omega$ consists of infinitely many connected components:}\ 
Note that $\omega \in \Omega_{\rm inf}$ in this case.
Now by Lemma \ref{lem:charofWtopm1}, we have $|{\rm Hom}(W_\omega,\{\pm 1\})|=\infty$ and so $|{\rm Hom}(W_{\rm inf},Z(W))|=\infty$ (since we assumed that $Z(W) \neq 1$).
Thus by (ii), the subgroup $H'_2H'_3H'_4$ ($\supset H$) has index $|H'_1|=\infty$, so that $H$ also has infinite index in ${\rm Aut}(W)$.\\ \indent
{\bf Subcase 2-2. The odd Coxeter graph of every $W_\omega$ consists of only finitely many connected components:}\ 
Then we have $|\Omega|=\infty$ by the hypothesis of Case 2.
Since we assumed that $Z(W) \neq 1$, we can take an infinite sequence $\omega_0,\omega_1,\omega_2,\dotsc$ of distinct elements of $\Omega$ such that $Z(W_{\omega_0}) \neq 1$.
Let $u$ denote the unique element of $Z(W_{\omega_0}) \smallsetminus \{1\}$.
Now for $k \geq 1$, we define $f_k \in {\rm Hom}(W,Z(W))$ componentwise by
\begin{equation*}
f_k(w)=
\begin{cases}
u^{\ell(w)} & \text{ if } \omega \in \{\omega_1,\dotsc ,\omega_k\} \text{ and } w \in W_\omega,\\
1 & \text{ if } \omega \in \Omega \smallsetminus \{\omega_1,\dotsc ,\omega_k\} \text{ and } w \in W_\omega.
\end{cases}
\end{equation*}
Then we have $f_k \circ f_k=1$ and so $f_k*f_k=1$ since $Z(W)$ is an elementary abelian $2$-group.
This implies that $f_k \in {\rm Hom}(W,Z(W))^\times$ and so ${f_k}^\flat \in {\rm Aut}(W)$, while ${\rm sep}({f_k}^\flat)=\{\omega_1,\dotsc ,\omega_k\}$ by definition.
Thus by the above remark, all ${f_k}^\flat$ belong to distinct cosets in ${\rm Aut}(W)/H$ and so $H$ has infinite index in ${\rm Aut}(W)$.
Hence the proof is concluded.
\end{proof}
\begin{exmp}
Let $m=(m_1,m_2,\dotsc)$ be an infinite sequence of nonnegative integers.
Here we examine ${\rm Aut}(W_m)$ for the group $W_m=\prod_{n \geq 1}({\rm Sym}_n)^{m_n}$ by using our result, where ${\rm Sym}_n={\rm Sym}(\{1,2,\dotsc ,n\})$ is the symmetric group of degree $n$.
Note that ${\rm Sym}_1=1$.\\ \indent
Since ${\rm Sym}_n$ ($n \geq 2$) is the Coxeter group $W(A_{n-1})$, which is directly indecomposable (cf.\ Theorem \ref{thm:indecomposabilityofW}), we can apply Theorem \ref{thm:automorphismgroup} (i) to this decomposition of $W_m$.
In this case, we have $Z({\rm Sym}_n)=1$ unless $Z({\rm Sym}_n)={\rm Sym}_n$ (namely $n=1,2$), so that ${\rm Hom}(W_m,Z(W_m))_o=1$.
Thus we have ${\rm Aut}(W_m)=H_1 \rtimes H_2 \rtimes H_3$.\\ \indent
Note that $Z(W_m)=({\rm Sym}_2)^{m_2} \simeq \{\pm 1\}^{m_2}$, while $|{\rm Hom}({\rm Sym}_n,\{\pm 1\})|=2$ for all $n \geq 2$ by Lemma \ref{lem:charofWtopm1}.
Thus Lemmas \ref{lem:Homofabelian} and \ref{lem:subgroupofHom} (ii) imply that
\begin{equation*}
\begin{split}
H_1&={\rm Hom}\biggl(\prod_{n \geq 3}({\rm Sym}_n)^{m_n},Z(W_m)\biggr)^\flat \rtimes {{\rm Hom}({\rm Sym}_2^{m_2},Z(W_m))^\times}^\flat\\
&=\left(\overline{\prod}_{n \geq 3}{\rm Hom}\bigl(({\rm Sym}_n)^{m_n},Z(W_m)\bigr)\right)^\flat \rtimes {\rm Aut}(({\rm Sym}_2)^{m_2})\\
& \simeq \left(\overline{\prod}_{n \geq 3}\{\pm 1\}^{m_2m_n}\right) \rtimes {\rm GL}_{m_2}(\mathbb{F}_2).
\end{split}
\end{equation*}
Secondly, recall the well-known fact that ${\rm Aut}({\rm Sym}_n)={\rm Inn}({\rm Sym}_n)$ (the group of inner automorphisms) if $n \neq 6$ and $|{\rm Aut}({\rm Sym}_6)/{\rm Inn}({\rm Sym}_6)|=2$.
This implies that ${\rm Aut}({\rm Sym}_2)=1$, $|{\rm Aut}({\rm Sym}_6)|=2|{\rm Sym}_6|$ and ${\rm Aut}({\rm Sym}_n) \simeq {\rm Sym}_n$ if $n \neq 2,6$.
Thus we have
\begin{equation*}
H_2 \simeq \overline{\prod}_{n \geq 3}{\rm Aut}({\rm Sym}_n)^{m_n} \simeq \left(\overline{\prod}_{3 \leq n \neq 6}{{\rm Sym}_n}^{m_n}\right) \times {\rm Aut}({\rm Sym}_6)^{m_6}.
\end{equation*}
Moreover, by definition, we have $H_3 \simeq \overline{\prod}_{n \geq 3}{\rm Sym}_{m_n}$.\\ \indent
As a special case, if all but finitely many terms in $m$ are $0$, then (by putting $|m|=\sum_n m_n<\infty$) we have
\begin{equation*}
\begin{split}
|H_1|&=2^{m_2(|m|-m_1-m_2)}\prod_{i=0}^{m_2-1}(2^{m_2}-2^i)=2^{m_2(|m|-m_1-m_2)+\binom{m_2}{2}}\prod_{i=1}^{m_2}(2^i-1),\\
|H_2|&=2^{m_6}\prod_{n \geq 3}(n!)^{m_n},\ |H_3|=\prod_{n \geq 3}m_n!.
\end{split}
\end{equation*}
Hence we have
\begin{equation*}
\begin{split}
|{\rm Aut}(W_m)|&=|H_1| \cdot |H_2| \cdot |H_3|\\
&=2^{m_2(|m|-m_1-m_2)+\binom{m_2}{2}+m_6}\prod_{i=1}^{m_2}(2^i-1)\prod_{n \geq 3}\left((n!)^{m_n}m_n!\right)\\
&=\biggl(2^{m_2(|m|-m_1-m_2-1)+\binom{m_2}{2}+m_6}\prod_{i=1}^{m_2}(2^i-1)\prod_{n \geq 3}m_n!\biggr)|W_m|.
\end{split}
\end{equation*}
\end{exmp}
%
\section{Centralizers of normal subgroups generated by involutions}\label{sec:ZofinvolutiveH}
\subsection{Proof of Theorem \ref{thm:Zofinvolutivenormalsubgroup}}\label{sec:proofoftheorem_Sec3}
In this section, we prove Theorem \ref{thm:Zofinvolutivenormalsubgroup}.
From now, $(W,S)$ always denotes a Coxeter system.
In the proof, we use the notion of core subgroups (cf.\ Section \ref{sec:generalgroups}).
For a subgroup $G \leq W$, let $X_G$ be the set of all elements in $G$ of the form $w_0(I)$ ($I \subset S$) such that $1 \neq w_0(I) \in Z(W_I)$.
Then we have the following relation (proved below):
\begin{prop}
\label{prop:ZofnormalsubgroupinW}
Let $H \lhd W$ be a normal subgroup generated by involutions.
Then $H$ is the smallest normal subgroup of $W$ containing $X_H$, and
\begin{equation*}
Z_W(H)=\bigcap_{w_0(I) \in X_H}{\rm Core}_W(N_W(W_I)).
\end{equation*}
\end{prop}
On the other hand, the subgroups ${\rm Core}_W(N_W(W_I))$ are determined completely (for irreducible $(W,S)$) by the following theorem, which we prove in later subsections.
Here we use the notation $(W(D_3),S(D_3))$ instead of $(W(A_3),S(A_3))$.
\begin{thm}
\label{thm:coreofnormalizer}
(cf.\ Definitions \ref{defn:coresubgroup} and \ref{defn:Coxetergraphs}.)
Let $(W,S)$ be an irreducible Coxeter system of an arbitrary rank, and $I$ nonempty proper subset of $S$.
Then:\\
{\bf (i)} If $(W,S)=(W(B_n),S(B_n))$, $1 \leq k<n \leq \infty$, $\tau \in {\rm Aut}(\varGamma(B_n))$ and $I=\tau(S(B_k))$, then ${\rm Core}_W(N_W(W_I))=\tau(G_{B_n})$.\\
{\bf (ii)} If $(W,S)=(W(D_n),S(D_n))$, $2 \leq k<n \leq \infty$, $\tau \in {\rm Aut}(\varGamma(D_n))$ and $I=\tau(S(D_k))$, then ${\rm Core}_W(N_W(W_I))=\tau(G_{D_n})$.\\
{\bf (iii)} Otherwise, ${\rm Core}_W(N_W(W_I))=Z(W)$.\\
(cf.\ Lemma \ref{lem:defofpropersubgroupG} for definition of $G_{B_n}$ and $G_{D_n}$.)
\end{thm}
Note that, on the other hand, ${\rm Core}_W(N_W(W_I))=N_W(W_I)=W$ if $I=\emptyset$ or $S$.
Theorem \ref{thm:Zofinvolutivenormalsubgroup} will be proved by combining Proposition \ref{prop:ZofnormalsubgroupinW} and Theorem \ref{thm:coreofnormalizer}.\\ \indent
In the proof of Proposition \ref{prop:ZofnormalsubgroupinW}, we use the following two results:
\begin{thm}[\cite{Ric}, Theorem A]
\label{thm:Richardson_involutions}
Let $w$ be an involution in $W$.
Then $w$ is conjugate in $W$ to some element $w_0(I)$ ($I \subset S$) such that $w_0(I) \in Z(W_I)$.
\end{thm}
\begin{lem}
\label{lem:Zoflongestelement}
Let $W_I$ be a finite parabolic subgroup of $W$ such that $w_0(I) \in Z(W_I)$.
Then $Z_W(w_0(I))=N_W(W_I)$.
\end{lem}
\begin{proof}
First, assume $u \in Z_W(w_0(I))$.
Then $u^{-1}w_0(I)u=w_0(I) \in Z(W_I)$ and so $w_0(I) \cdot (u \cdot \alpha_s)=uw_0(I) \cdot \alpha_s=-u \cdot \alpha_s$ for all $s \in I$.
This implies that $u \cdot \alpha_s \in \Phi_I$ for all $s \in I$, so that $u \in N_W(W_I)$ by (\ref{eq:charofnormalizer}).\\ \indent
Conversely, assume $u \in N_W(W_I)$.
Put $u'=uw_0(I)u^{-1} \in W_I$.
Then we have $u' \cdot \alpha_s=-\alpha_s$ for all $s \in I$ (since $w_0(I)$ maps $u^{-1} \cdot \alpha_s \in \Phi_I$ (cf.\ (\ref{eq:charofnormalizer})) to $-u^{-1} \cdot \alpha_s$).
Hence we have $u'=w_0(I)$ and so $u \in Z_W(w_0(I))$.
\end{proof}
\begin{proof}[Proof of Proposition \ref{prop:ZofnormalsubgroupinW}]
By Theorem \ref{thm:Richardson_involutions}, every involution in $H$ is conjugate to some element of $X_H$ (since $H \lhd W$).
This implies that any normal subgroup of $W$ containing $X_H$ also contains all the generators of $H$.
Thus the first claim follows.
For the second one, apply Lemmas \ref{lem:Zandcore} and \ref{lem:Zoflongestelement}.
\end{proof}
\begin{proof}[Proof of Theorem \ref{thm:Zofinvolutivenormalsubgroup}]
The claim (i) is obvious.
From now, we assume $H \not\subset Z(W)$.
Note that $Z(W) \subset Z_W(H)$.
Note also that, by Proposition \ref{prop:ZofnormalsubgroupinW},
\begin{equation}
\label{eq:keyofthm_Zofinvolutivenormalsubgroup}
Z_W(H) \subset {\rm Core}_W(N_W(W_I)) \text{ for all } w_0(I) \in X_H.
\end{equation}
\indent
{\bf Case 1. $(W,S)=(W(B_n),S(B_n))$, $n \geq 2$ or $(W(D_n),S(D_n))$, $n \geq 3$:}\ 
Let $\mathcal{T}=B$, $L=1$ for the former case, $\mathcal{T}=D$, $L=2$ for the latter case.\\ \indent
{\bf Subcase 1-1. $\mathcal{T}=B$, $n \neq 2$ or $\mathcal{T}=D$, $n \neq 4$:}\ 
Note that in this case, any automorphism of $\varGamma(\mathcal{T}_n)$ preserves the sets $S(\mathcal{T}_k)$, elements $w_0(S(\mathcal{T}_k))$ ($k \geq L$) and so the subgroup $G_{\mathcal{T}_n}$.\\ \indent
{\bf Subsubcase 1-1-1. $H \subset G_{\mathcal{T}_n}$:}\ 
This is a case (ii) or (iii) (for $\tau$ identity), and so we have to show $Z_W(H)=G_{\mathcal{T}_n}$.
The inclusion $\supset$ holds since $G_{\mathcal{T}_n}$ is abelian.
Conversely, since $H \not\subset Z(W)$, $X_H$ contains an element other than $w_0(S)$, so that we have $Z_W(H) \subset G_{\mathcal{T}_n}$ by (\ref{eq:keyofthm_Zofinvolutivenormalsubgroup}) and Theorem \ref{thm:coreofnormalizer}.\\ \indent
{\bf Subsubcase 1-1-2. $H \not\subset G_{\mathcal{T}_n}$:}\ 
By the above remark, this is actually not a case (ii) or (iii), so that we have to show $Z_W(H) \subset Z(W)$.
Now $X_H$ contains an element $w_0(I)$ such that $I \neq S(\mathcal{T}_k)$ for any $L \leq k \leq n$, since otherwise $H \subset G_{\mathcal{T}_n}$ by Lemma \ref{lem:defofpropersubgroupG}.
For this $I$, we have ${\rm Core}_W(N_W(W_I))=Z(W)$ by Theorem \ref{thm:coreofnormalizer}, so that the claim follows from (\ref{eq:keyofthm_Zofinvolutivenormalsubgroup}).\\ \indent
{\bf Subcase 1-2. $\mathcal{T}=B$, $n=2$:}\ 
Note that $X_H \subset \{s_1,s_2,w_0(S)\}$ in this case.
Moreover, $X_H \not\subset \{w_0(S)\}$ since $H \not\subset Z(W)$.\\ \indent
{\bf Subsubcase 1-2-1. $s_1 \in X_H$ and $s_2 \not\in X_H$:}\ 
In this case, we have $X_H \subset \{s_1,w_0(S)\}$ and so $H \subset G_{B_2}$ by Lemma \ref{lem:defofpropersubgroupG}.
This is a case (ii) (for $\tau$ identity).
Now we have $G_{B_2} \subset Z_W(H)$ since $G_{B_2}$ is abelian, while $Z_W(H) \subset G_{B_2}$ by (\ref{eq:keyofthm_Zofinvolutivenormalsubgroup}) and Theorem \ref{thm:coreofnormalizer} (applying to $\{s_1\} \subset S$).
Thus the claim holds.\\ \indent
{\bf Subsubcase 1-2-2. $s_1 \not\in X_H$ and $s_2 \in X_H$:}\ 
By symmetry, this is also a case (ii) (for the unique $\tau \neq {\rm id}_S$) and the claim holds similarly.\\ \indent
{\bf Subsubcase 1-2-3. $s_1 \in X_H$ and $s_2 \in X_H$:}\ 
Note that $H=W$.
This is not a case (ii) or (iii), and actually $Z_W(H)=Z(W)$.\\ \indent
{\bf Subcase 1-3. $\mathcal{T}=D$, $n=4$:}\ 
Note that (by definition)
\begin{equation*}
X_H \subset \{s_1,\ s_2,\ s_3,\ s_4,\ s_1s_2s_4,\ s_1s_2,\ s_2s_4,\ s_4s_1,\ w_0(S)\}.
\end{equation*}
\indent
{\bf Subsubcase 1-3-1. $X_H$ contains one of the first five elements:}\ 
Now we have $H \not\subset \tau(G_{D_4})$ for any $\tau$, so that this is not a case (iii) and we have to show $Z_W(H) \subset Z(W)$.
This claim follows from (\ref{eq:keyofthm_Zofinvolutivenormalsubgroup}) (applying to the element of $X_H$ given in the hypothesis here) and Theorem \ref{thm:coreofnormalizer}.\\ \indent
{\bf Subsubcase 1-3-2. $X_H$ contains at least two of the elements $s_1s_2$, $s_2s_4$, $s_4s_1$:}\ 
Now we have $H \not\subset \tau(G_{D_4})$ for any $\tau$, so that this is not a case (iii) and we have to show $Z_W(H) \subset Z(W)$.
Let $X_H$ contain two such elements $s_is_j$, $s_js_k$, and put $I=\{s_i,s_j\}$, $J=\{s_j,s_k\}$.
Then we have
\begin{equation*}
{\rm Core}_W(N_W(W_I)) \cap {\rm Core}_W(N_W(W_J)) \subset {\rm Core}_W(N_W(W_{\{s_j\}}))
\end{equation*}
by (\ref{eq:coreofintersection}), (\ref{eq:intersectionofnormalizer}) and (\ref{eq:coreisincreasing}).
Thus we have $Z_W(H) \subset {\rm Core}_W(N_W(W_{\{s_j\}}))=Z(W)$ by (\ref{eq:keyofthm_Zofinvolutivenormalsubgroup}) and Theorem \ref{thm:coreofnormalizer}.\\ \indent
{\bf Subsubcase 1-3-3. $X_H$ contains none of the first five elements and at most one of $s_1s_2$, $s_2s_4$, $s_4s_1$:}\ 
Note that $X_H \not\subset \{w_0(S)\}$ since $H \not\subset Z(W)$.
Thus we have $s_is_j \in X_H \subset \{s_is_j,w_0(S)\}$ for one of $(i,j)=(1,2)$, $(2,4)$, $(4,1)$.
Lemma \ref{lem:defofpropersubgroupG} implies that this is a case (iii) (namely $H \subset \tau(G_{D_4})$), by taking $\tau \in {\rm Aut}(\varGamma)$ mapping $s_1$, $s_2$ to $s_i$, $s_j$ respectively.
Now $\tau(G_{D_4}) \subset Z_W(H)$ since $\tau(G_{D_4})$ is abelian.
Conversely, we have ${\rm Core}_W(N_W(W_{\{s_i,s_j\}}))=\tau(G_{D_4})$ by Theorem \ref{thm:coreofnormalizer}, so that $Z_W(H) \subset \tau(G_{D_4})$ by (\ref{eq:keyofthm_Zofinvolutivenormalsubgroup}).
Thus the claim holds.\\ \indent
{\bf Case 2. $(W,S) \not\simeq (W(B_n),S(B_n))$ ($n \geq 2$), $(W(D_n),S(D_n))$ ($n \geq 3$):}\ 
This is not a case (ii) or (iii), so that we have to show $Z_W(H) \subset Z(W)$.
Since $H \not\subset Z(W)$, $X_H$ contains an element other than $w_0(S)$, so that we have $Z_W(H) \subset Z(W)$ by (\ref{eq:keyofthm_Zofinvolutivenormalsubgroup}) and Theorem \ref{thm:coreofnormalizer}.
Hence the proof is concluded.
\end{proof}
%
\subsection{Some lemmas}\label{sec:lemmaforcores}
In the rest of this paper, we prove Theorem \ref{thm:coreofnormalizer}.
In this subsection, we prepare some lemmas used in our proof.
From now, we abbreviate the notation ${\rm Core}_W(N_W(W_I))$ to $C_I$.\\ \indent
First, by combining Lemma \ref{lem:intersectionofnormalizers}, (\ref{eq:coreofintersection}) and (\ref{eq:coreisincreasing}), we have:
\begin{equation}
\label{eq:coreofseparatedI}
\text{ If } I \subset J \subset S \text{ and } J \smallsetminus I \subset I^\perp, \text{ then } C_J \cap C_I \subset C_{J \smallsetminus I}.
\end{equation}
\begin{lem}[Expanding Lemma]
If $I \subset S$ and $s \in S \smallsetminus (I \cup I^\perp)$, then $C_I \subset C_{I \cup \{s\}}$.
\end{lem}
\begin{proof}
It is enough (by (\ref{eq:coreinsubgroup})) to show that $C_I \subset N_W(W_{I \cup \{s\}})$.
Let $w \in C_I$.
By the hypothesis, we have $c=\langle \alpha_s,\alpha_t \rangle <0$ for some $t \in I$.
Now since $sws \in C_I \subset N_W(W_I)$, we have $sws \cdot \alpha_t \in \Phi_I$ (by (\ref{eq:charofnormalizer})) and so $ws \cdot \alpha_t \in \Phi_{I \cup \{s\}}$.
On the other hand, we have $ws \cdot \alpha_t=w \cdot \alpha_t-2cw \cdot \alpha_s$.
Thus $w \cdot \alpha_s \in \Phi_{I \cup \{s\}}$ since $w \cdot \alpha_t \in \Phi_I$ (by (\ref{eq:charofnormalizer})).
Hence we have $w \in N_W(W_{I \cup \{s\}})$ by (\ref{eq:charofnormalizer}).
\end{proof}
For $s \in S$ and $I \subset S$, let $d_{\varGamma}(s,I)=\min\{d_{\varGamma}(s,t) \mid t \in I\}$ denote the distance from $s$ to the set $I$ in the Coxeter graph $\varGamma$ of $(W,S)$.
\begin{lem}[Cutting Lemma]
Let $(W,S)$ be irreducible, $I \subset S$ and $s \in S \smallsetminus I$.
Then for $d_{\varGamma}(s,I)<k<\infty$, we have $C_I \subset C_J$, where $J=\{t \in I \mid d_{\varGamma}(s,t) \geq k\}$.
\end{lem}
\begin{proof}
It is enough (by (\ref{eq:coreinsubgroup}) and (\ref{eq:charofnormalizer})) to show that $w \cdot \Phi_J \subset \Phi_J$ (or equivalently, $w \cdot \Pi_J \subset \Phi_J$) for all $w \in C_I$.
Assume contrary that $t \in J$ and $w \cdot \alpha_t \not\in \Phi_J$.
Note that $w \cdot \alpha_t \in \Phi_I$ (by (\ref{eq:charofnormalizer})) and so $s \not\in {\rm supp}(w \cdot \alpha_t)$.
Then by definition of $J$, we have
\begin{equation*}
\text{($d=$) } d_{\varGamma}(s,{\rm supp}(w \cdot \alpha_t))<k \leq d_{\varGamma}(s,t).
\end{equation*}
Take a shortest path $s_0=s,s_1,\dotsc ,s_{d-1},s_d \in {\rm supp}(w \cdot \alpha_t)$ in $\varGamma$ from $s$ to the set ${\rm supp}(w \cdot \alpha_t)$.
Then by the above inequality, we have $s_i \in \{t\}^\perp$ for all $0 \leq i \leq d-1$.
Put $u=ss_1 \dotsm s_{d-1} \in W$.
Then we have $uwu^{-1} \cdot \alpha_t=uw \cdot \alpha_t$ and so (by (\ref{eq:rootandadjacentelement}))
\begin{equation*}
{\rm supp}(uwu^{-1} \cdot \alpha_t)={\rm supp}(w \cdot \alpha_t) \cup \{s,s_1,\dotsc ,s_{d-1}\} \not\subset I
\end{equation*}
(note that $s \not\in I$).
On the other hand, we have $uwu^{-1} \in C_I$ and so $uwu^{-1} \cdot \alpha_t \in \Phi_I$ (by (\ref{eq:charofnormalizer})).
This is a contradiction.
Hence the claim holds.
\end{proof}
\begin{lem}[Shifting Lemma]
Suppose that $s,t \in S$ are in the same connected component of the odd Coxeter graph $\varGamma^{\rm odd}$ of $(W,S)$.
Then $C_{\{s\}}=C_{\{t\}}$.
\end{lem}
\begin{proof}
By definition of $\varGamma^{\rm odd}$, and by symmetry, it is enough to show that $C_{\{s\}} \subset C_{\{t\}}$ for any $s,t$ such that $m(s,t)=2k+1$ is odd.
Now by putting $u=(st)^k \in W$, we have $t=usu^{-1}$.
Thus for $w \in C_{\{s\}}$, we have
\begin{equation*}
wtw^{-1}=wusu^{-1}w^{-1}=u(u^{-1}wu)s(u^{-1}wu)^{-1}u^{-1}=usu^{-1}=t
\end{equation*}
since $u^{-1}wu \in C_{\{s\}}$.
Thus $w \in N_W(W_{\{t\}})$.
Hence the claim follows from (\ref{eq:coreinsubgroup}).
\end{proof}
Moreover, we have:
\begin{lem}
\label{lem:coreofparabolic}
Let $(W,S)$ be irreducible and $I$ a nontrivial proper subset of $S$.
Then ${\rm Core}_W(W_I)=1$.
\end{lem}
\begin{proof}
Assume contrary that $1 \neq w \in {\rm Core}_W(W_I)$ (so that $w \cdot \Phi_I=\Phi_I$ by (\ref{eq:charofnormalizer})).
Fix $s \in S \smallsetminus I$ and take $\gamma \in \Phi_I^+$ such that $w \cdot \gamma \in \Phi_I^-$.\\ \indent
{\bf Case 1. ($d=$) $d_{\varGamma}(s,{\rm supp}(\gamma)) \leq d_{\varGamma}(s,{\rm supp}(w \cdot \gamma))$:}\ 
Take a shortest path $s_0=s,s_1,\dotsc ,s_{d-1},s_d \in {\rm supp}(\gamma)$ in $\varGamma$ from $s$ to the set ${\rm supp}(\gamma)$.
Then by the above inequality, we have $s_i \not\in {\rm supp}(w \cdot \gamma)$ for all $0 \leq i \leq d-1$.
Put $u=ss_1 \dotsm s_{d-1} \in W$.
Then we have $u \cdot \gamma \in \Phi^+$ (by (\ref{eq:rootandsmallerparabolic})), ${\rm supp}(u \cdot \gamma)={\rm supp}(\gamma) \cup \{s,s_1,\dotsc ,s_{d-1}\} \not\subset I$ (by (\ref{eq:rootandadjacentelement})) and so $u \cdot \gamma \in \Phi^+ \smallsetminus \Phi_I$.
On the other hand, we have $uwu^{-1} \cdot (u \cdot \gamma)=u \cdot (w \cdot \gamma) \in \Phi^-$ (by (\ref{eq:rootandsmallerparabolic})).
This is a contradiction, since $uwu^{-1} \in {\rm Core}_W(W_I) \subset W_I$.\\ \indent
{\bf Case 2. $d_{\varGamma}(s,{\rm supp}(\gamma))>d_{\varGamma}(s,{\rm supp}(w \cdot \gamma))$:}\ 
Now by applying Case 1 to the elements $w^{-1} \in {\rm Core}_W(W_I)$ and $-w \cdot \gamma \in \Phi_I\left[\!\right.w^{-1}\left.\!\right]$, we have a contradiction again.
Hence the claim holds in any case.
\end{proof}
Owing to Lemma \ref{lem:coreofparabolic}, we have the following results:
\begin{align}
\label{eq:coreofmaximalparabolic}
&\text{ If } (W,S) \text{ is irreducible, } |W|=\infty \text{ and } s \in S, \text{ then } C_{S \smallsetminus \{s\}}=1.\\
\label{eq:coreofinfinite}
&\text{ If } I \text{ is an irreducible component of } J \subset S \text{ and } |W_I|=\infty, \text{ then } C_J=1.
\end{align}
(Here we use Corollary \ref{cor:Nofmaximalparabolic} (iii), Proposition \ref{prop:Nofinfiniteirreducible}, respectively.)
\subsection{Proof for finite case}\label{sec:coreoffinitecase}
In this subsection, we prove Theorem \ref{thm:coreofnormalizer} for the case $|W|<\infty$.
From now, we abbreviate often the terms ``Expanding Lemma'', ``Cutting Lemma'', ``Shifting Lemma'' to `EL', `CL', `SL', respectively.
\begin{lem}
\label{lem:coreofZofsinglereflection}
Let $(W,S)$ be irreducible, $|W|<\infty$ and $s \in S$.
Suppose that no condition below is satisfied: {\bf (I)} $W=W(B_n)$, $n \geq 2$, $s=s_1$, {\bf (II)} $W=W(B_2)$, $s=s_2$, {\bf (III)} $W=W(I_2(m))$, $m$ even.
Then $C_{\{s\}}=Z(W)$.
\end{lem}
\begin{proof}
Since $Z(W) \subset C_{\{s\}}$ and $\bigcap_{t \in S}N_W(W_{\{t\}})=Z(W)$, it is enough to show that $C_{\{s\}} \subset C_{\{t\}}$ for all $t \in S$.\\ \indent
{\bf Case 1. The odd Coxeter graph $\varGamma^{\rm odd}$ of $(W,S)$ is connected:}\ 
Then the claim follows from the Shifting Lemma.\\ \indent
{\bf Case 2. $W=W(B_n)$, $n \geq 3$ and $s \neq s_1$:}\ 
We have $C_{\{s\}} \overset{\text{SL}}{=} C_{\{s_i\}}$ for all $2 \leq i \leq n$, while $C_{\{s_2\}} \overset{\text{EL}}{\subset} C_{\{s_1,s_2\}} \overset{\text{CL}}{\subset} C_{\{s_1\}}$ (since $n \geq 3$).
Thus the claim holds.\\ \indent
{\bf Case 3. $W=W(F_4)$:}\ 
By symmetry, we may assume $s=s_1$ or $s_2$.
Now we have $C_{\{s_1\}} \overset{\text{SL}}{=} C_{\{s_2\}} \overset{\text{EL}}{\subset} C_{\{s_2,s_3\}} \overset{\text{CL}}{\subset} C_{\{s_3\}} \overset{\text{SL}}{=} C_{\{s_4\}}$.
Hence the claim holds.
\end{proof}
\begin{cor}
\label{cor:farthestvertex}
Let $(W,S)$ be irreducible, $|W|<\infty$, $s \in S$ and suppose that there is a unique vertex $t$ of $\varGamma$ farthest from $s$.
Suppose further that $W$ and $t$ do not satisfy any of the three conditions (I)--(III) in Lemma \ref{lem:coreofZofsinglereflection}.
Then $C_{S \smallsetminus \{s\}}=Z(W)$.
\end{cor}
\begin{proof}
Now we have $C_{S \smallsetminus \{s\}} \overset{\text{CL}}{\subset} C_{\{t\}}$ by the choice of $t$.
Then apply Lemma \ref{lem:coreofZofsinglereflection}.
\end{proof}
\begin{lem}
\label{lem:coreforsmallS}
Suppose that one of the following conditions is satisfied: {\bf (I)} $W=W(B_3)$, $s=s_2$, {\bf (II)} $W=W(D_4)$, $s=s_3$, {\bf (III)} $W=W(H_3)$, $s=s_2$, {\bf (IV)} $W=W(I_2(m))$ ($m \geq 6$ even), $s \in S$.
Then $C_I=Z(W)$, where $I=S \smallsetminus \{s\}$.
\end{lem}
\begin{proof}
By the hypothesis and Corollary \ref{cor:Nofmaximalparabolic} (ii), we have $N_W(W_I)=W_I \times Z(W)$.
Now a direct computation shows that $sW_Is \cap N_W(W_I)=1$, so that $W_I \cap C_I=1$ by (\ref{eq:subgroupoutsidecore}).
Since $Z(W) \subset C_I$, we have $C_I=Z(W)$.
\end{proof}
\begin{lem}
\label{lem:coreofspecialsubgroup}
{\bf (i)} If $W=W(B_n)$, $1 \leq n<\infty$, then ${\rm Core}_W(G_{B_n})=G_{B_n}$.\\
{\bf (ii)} If $W=W(D_n)$, $3 \leq n<\infty$, then ${\rm Core}_W(G_{D_n} \rtimes \langle s_1 \rangle)=G_{D_n}$.
\end{lem}
\begin{proof}
The claim (i) is obvious, since $G_{B_n} \lhd W$ (cf.\ Lemma \ref{lem:defofpropersubgroupG}).
For (ii), we have $G_{D_n} \subset {\rm Core}_W(G_{D_n} \rtimes \langle s_1 \rangle)$ since $G_{D_n} \lhd W$, while $s_1 \not\in {\rm Core}_W(G_{D_n} \rtimes \langle s_1 \rangle)$ since $s_1s_3s_1s_3s_1=s_3 \not\in G_{D_n} \rtimes \langle s_1 \rangle$.
Thus the claim holds.
\end{proof}
\begin{proof}[Proof of Theorem \ref{thm:coreofnormalizer} (for finite $W$)]
Note that $Z(W) \subset C_I$ by definition.\\ \indent
{\bf Case 1. $(W,S)=(W(\mathcal{T}_n),S(\mathcal{T}_n))$ for $\mathcal{T}=B$, $n \geq 3$ or $\mathcal{T}=D$, $3 \leq n \neq 4$:}\ 
Put $L=1$ in the former case, $L=2$ in the latter case.
Note that in this case, any automorphism of $\varGamma(\mathcal{T}_n)$ preserves the sets $S(\mathcal{T}_k)$, elements $w_0(S(\mathcal{T}_k))$ ($k \geq L$) and so the subgroup $G_{\mathcal{T}_n}$.\\ \indent
{\bf Subcase 1-1. $I=S(\mathcal{T}_k)$ for some $L \leq k<n$:}\ 
This is a case (i) or (ii) of Theorem \ref{thm:coreofnormalizer} (for $\tau$ identity), so that we have to show $C_I=G_{\mathcal{T}_n}$.
Note that
\begin{equation*}
C_{S(\mathcal{T}_i)} \overset{\text{EL}}{\subset} C_{S(\mathcal{T}_j)} \overset{\text{CL}}{\subset} C_{S(\mathcal{T}_i)} \text{ and so } C_{S(\mathcal{T}_i)}=C_{S(\mathcal{T}_j)} \text{ for all } L \leq i<j<n.
\end{equation*}
Thus we may assume $I=S(\mathcal{T}_L)$, and we have $C_I \subset \bigcap_{i=L}^{n-1}N_W(W_{S(\mathcal{T}_i)})$.
By Corollary \ref{cor:Nforspecialcase}, (\ref{eq:coreinsubgroup}) and Lemma \ref{lem:coreofspecialsubgroup}, we have $C_I \subset G_{\mathcal{T}_n}$.
Conversely, since $G_{\mathcal{T}_n}$ is abelian and contains $w_0(I)$, we have $G_{\mathcal{T}_n} \subset Z_W(w_0(I))=N_W(W_I)$ by Lemma \ref{lem:Zoflongestelement}.
Thus $G_{\mathcal{T}_n} \subset C_I$ since $G_{\mathcal{T}_n} \lhd W$.
Hence $C_I=G_{\mathcal{T}_n}$.\\ \indent
{\bf Subcase 1-2. $I \neq S(\mathcal{T}_k)$ for all $L \leq k<n$:}\ 
By the above remark, this is not a case (i) or (ii), and so we have to show $C_I \subset Z(W)$.
Note that $I \neq S$.
Let $M$ be the first index $\geq 1$ such that $s_M \not\in I$, so that $S(\mathcal{T}_{M-1}) \subset I$ (where we put $S(\mathcal{T}_0)=\emptyset$).
If $\mathcal{T}=D$ and $M=2$, then we have $C_I \overset{\text{EL}}{\subset} C_{S \smallsetminus \{s_M\}}$ since $I \neq \emptyset$.
Otherwise, there is some $M<i \leq n$ such that $s_i \in I$ (since otherwise we have a contradiction $I=S(\mathcal{T}_{M-1})$), and so $M<n$ and $C_I \overset{\text{EL}}{\subset} C_{S \smallsetminus \{s_M\}}$.
In any case, we may assume that $I=S \smallsetminus \{s_M\}$.
Now there are the following three cases:\\ \indent
{\bf Subsubcase 1-2-1. $M \leq L+1$:}\ 
Note that $M<n$, and so $(\mathcal{T}_n,M) \neq (D_3,3)$.
If $\mathcal{T}_n=B_3$ and $M=2$, then $C_I=Z(W)$ by Lemma \ref{lem:coreforsmallS}.
Otherwise, we have a unique vertex of $\varGamma$ farthest from $s$; that is $s_{3-M}$ if $\mathcal{T}_n=D_3$ and $M \leq 2$, and $s_n$ otherwise (note that $\mathcal{T}_n \neq D_4$).
Thus $C_I=Z(W)$ by Corollary \ref{cor:farthestvertex}.\\ \indent
{\bf Subsubcase 1-2-2. $L+2 \leq M \leq n-2$:}\ 
This hypothesis implies that
\begin{equation*}
C_I \overset{\text{CL}}{\subset} C_{I \smallsetminus \{s_{M-1},s_{M+1}\}} \overset{\text{EL}}{\subset} C_{S \smallsetminus \{s_{M-1}\}},
\end{equation*}
so that the claim follows inductively from the case of smaller $M$.\\ \indent
{\bf Subsubcase 1-2-3. $L+2 \leq M=n-1$:}\ 
Note that $n \geq L+3$ and $I=S(\mathcal{T}_{n-2}) \cup \{s_n\}$.
Now we have $C_I \overset{\text{CL}}{\subset} C_{S(\mathcal{T}_{n-3})} \overset{\text{EL}}{\subset} C_{S(\mathcal{T}_{n-2})}$ and so $C_I \subset C_{\{s_n\}}$ by (\ref{eq:coreofseparatedI}).
Thus $C_I \subset C_{\{s_n\}}=Z(W)$ by Lemma \ref{lem:coreofZofsinglereflection}.\\ \indent
{\bf Case 2. $(W,S)=(W(B_2),S(B_2))$:}\ 
Since $I$ is proper and nonempty, we have $I=\{s_i\}$ ($i=1$ or $2$).
This is a case (i), by taking $\tau={\rm id}_S$ (if $i=1$), $\tau \neq {\rm id}_S$ (if $i=2$).
Now we have to show $C_I=\tau(G_{B_2})$.
We have $C_I \subset N_W(W_{\tau(\{s_1\})})=\tau(G_{B_2})$ by Corollary \ref{cor:Nforspecialcase} (i).
Conversely, we have $\tau(G_{B_2}) \subset C_I$ by a similar argument to Subcase 1-1.
Thus $C_I=\tau(G_{B_2})$.\\ \indent
{\bf Case 3. $(W,S)=(W(D_4),S(D_4))$:}\ 
Note that $I$ is proper and nonempty.\\ \indent
{\bf Subcase 3-1. $|I|=1$:}\ 
This is not a case (i) or (ii), so that we have to show $C_I \subset Z(W)$.
This follows from Lemma \ref{lem:coreofZofsinglereflection}.\\ \indent
{\bf Subcase 3-2. $|I|=2$ and $s_3 \in I$:}\ 
This is also not a case (i) or (ii), so that we have to show $C_I \subset Z(W)$.
Let $I=\{s_3,s_i\}$.
Then we have $C_I \overset{\text{CL}}{\subset} C_{\{s_i\}}$, while $C_{\{s_i\}}=Z(W)$ by the previous case.
Thus $C_I \subset Z(W)$.\\ \indent
{\bf Subcase 3-3. $|I|=2$ and $s_3 \not\in I$:}\ 
Note that there is $\tau \in {\rm Aut}(\varGamma)$ such that $\tau(S(D_2))=I$.
This is a case (ii), so that we have to show $C_I=\tau(G_{D_4})$.
By symmetry, we may assume $\tau={\rm id}_S$.
First, we have $C_I \overset{\text{EL}}{\subset} C_{S(D_3)}$ and so $C_I \subset \bigcap_{i=2}^{3}N_W(W_{S(D_i)})=G_{D_4} \rtimes \langle s_1 \rangle$ by Corollary \ref{cor:Nforspecialcase} (ii).
Thus we have $C_I \subset G_{D_4}$ by (\ref{eq:coreinsubgroup}) and Lemma \ref{lem:coreofspecialsubgroup}.
Conversely, we have $G_{D_4} \subset C_I$ by a similar argument to Subcase 1-1.
Hence we have $C_I=G_{D_4}$.\\ \indent
{\bf Subcase 3-4. $|I|=3$ and $s_3 \in I$:}\ 
Note that there is $\tau \in {\rm Aut}(\varGamma)$ such that $\tau(S(D_3))=I$.
This is a case (ii), so that we have to show $C_I=\tau(G_{D_4})$.
By symmetry, we may assume $\tau={\rm id}_S$.
Now we have $C_I \overset{\text{CL}}{\subset} C_{S(D_2)} \overset{\text{EL}}{\subset} C_I$, while $C_{S(D_2)}=G_{D_4}$ by the previous subcase.
Thus $C_I=G_{D_4}$.\\ \indent
{\bf Subcase 3-5. $I=S \smallsetminus \{s_3\}$:}\ 
This is not a case (i) or (ii), so that we have to show $C_I \subset Z(W)$.
This follows from Lemma \ref{lem:coreforsmallS}.\\ \indent
{\bf Case 4. $(W,S) \not\simeq (W(B_n),S(B_n))$ ($n \geq 2$), $(W(D_n),S(D_n))$ ($n \geq 3$):}\ 
This is not a case (i) or (ii), so that we have to show $C_I \subset Z(W)$.
Note that $|S| \geq 2$ since $I$ is proper and nonempty.\\ \indent
{\bf Subcase 4-1. $|S|=2$:}\ 
Namely, $(W,S)=(W(\mathcal{T}),S(\mathcal{T}))$, $\mathcal{T}=A_2$ or $I_2(m)$ ($5 \leq m<\infty$), and $|I|=1$.
Then we have $C_I=Z(W)$ by Lemma \ref{lem:coreforsmallS} (for the latter case, with $m$ even) or Lemma \ref{lem:coreofZofsinglereflection} (the other cases).\\ \indent
{\bf Subcase 4-2. $|S|=3$:}\ 
Namely, $(W,S)=(W(H_3),S(H_3))$ (note that $W(A_3) \simeq W(D_3)$).
Now we have $C_I \overset{\text{EL}}{\subset} C_{S \smallsetminus \{s_i\}}$ for some $i$, while $C_{S \smallsetminus \{s_i\}}=Z(W)$ by Lemma \ref{lem:coreforsmallS} (if $i=2$) or Corollary \ref{cor:farthestvertex} (if $i \neq 2$).
Thus $C_I \subset Z(W)$.\\ \indent
{\bf Subcase 4-3. $|S| \geq 4$:}\ 
Namely, $(W,S)=(W(\mathcal{T}),S(\mathcal{T}))$ for $\mathcal{T}=A_n$ ($n \geq 4$), $E_n$ ($n=6,7,8$), $F_4$ or $H_4$.
Now we have $C_I \overset{\text{EL}}{\subset} C_{S \smallsetminus \{s_i\}}$ for some $i$.
Thus we may assume $I=S \smallsetminus \{s_i\}$.\\ \indent
{\bf Subsubcase 4-3-1. There is a unique vertex of $\varGamma$ farthest from $s_i$:}\ 
Now we have $C_I=Z(W)$ by Corollary \ref{cor:farthestvertex}.\\ \indent
{\bf Subsubcase 4-3-2. There are at least two vertices of $\varGamma$ farthest from $s_i$:}\ 
Namely, we have $(\mathcal{T},i)=(A_{2k+1},k+1)$ ($k \geq 2$), $(E_6,2)$, $(E_6,4)$ or $(E_8,5)$.
Now there are exactly two vertices $s,t$ of $\varGamma$ farthest from $s_i$, and there is a vertex $\neq s,t$ adjacent to $s$ and not adjacent to $t$.
This implies that $C_I \overset{\text{CL}}{\subset} C_{\{s,t\}} \overset{\text{CL}}{\subset} C_{\{t\}}$, while $C_{\{t\}}=Z(W)$ by Lemma \ref{lem:coreofZofsinglereflection}.
Thus $C_I \subset Z(W)$.
Hence the proof is concluded.
\end{proof}
%
\subsection{Proof for infinite case}\label{sec:coreofinfinitecase}
In this subsection, we prove Theorem \ref{thm:coreofnormalizer} in the case $|W|=\infty$.
The key facts are (\ref{eq:coreofmaximalparabolic}) and (\ref{eq:coreofinfinite}).\\ \indent
In the proof, we use a characterization (Proposition \ref{prop:nofinitecomponent}) of certain infinite Coxeter systems, which is based on the characterization of connected Coxeter graphs of finite type.
Before stating this, we prepare the following graph-theoretic lemma.
\begin{lem}
\label{lem:lemmafortree}
Let $\mathcal{G}$ be a connected acyclic graph (i.e.\ a tree) on nonempty vertex set $V(\mathcal{G})$ of an arbitrary cardinality (with no edge labels here).\\
{\bf (i)} If all vertices of $\mathcal{G}$ have degree $\leq 2$ and $\mathcal{G}$ has a terminal vertex (i.e.\ vertex of degree $1$) $s_0$, then $\mathcal{G} \simeq \varGamma(A_n)$ (as unlabelled graphs) for some $1 \leq n \leq \infty$.\\
{\bf (ii)} If $s_0 \in V(\mathcal{G})$ and all vertices of $\mathcal{G}$ except $s_0$ have degree $\leq 2$, then each connected component $\mathcal{G}'$ of $\mathcal{G} \smallsetminus \{s_0\}$ contains exactly one vertex $s$ adjacent to $s_0$, $\mathcal{G}' \simeq \varGamma(A_n)$ (as unlabelled graphs) for some $1 \leq n \leq \infty$ and $s$ is a terminal vertex of $\mathcal{G}'$.\\
{\bf (iii)} If all vertices of $\mathcal{G}$ have degree $2$, then $\mathcal{G} \simeq \varGamma(A_{\infty,\infty})$ (as unlabelled graphs).
\end{lem}
\begin{proof}
{\bf (i)} By the hypothesis, for any $s \in V(\mathcal{G})$, $\mathcal{G}$ contains a unique simple path $P_s=(t_s^{(0)}=s_0,t_s^{(1)},\dotsc ,t_s^{(\ell-1)},t_s^{(\ell)}=s)$ from $s_0$ to $s$.
Let $\ell(s)=\ell$, the length of $P_s$.
Then for all $s_1,s_2 \in V(\mathcal{G})$, we have either $P_{s_1} \subset P_{s_2}$ or $P_{s_2} \subset P_{s_1}$: Otherwise, for the first index $k$ such that $t_{s_1}^{(k)} \neq t_{s_2}^{(k)}$, the vertex $t_{s_1}^{(k-1)}=t_{s_2}^{(k-1)}$ is adjacent to distinct vertices $t_{s_1}^{(k)}$, $t_{s_2}^{(k)}$ (and $t_{s_1}^{(k-2)}$ if $k \geq 2$) but this is impossible by the hypothesis on the degree of $t_{s_1}^{(k-1)}$.\\ \indent
This observation shows that the map $\ell:V(\mathcal{G}) \to \{0,1,2,\dotsc\}$ is injective and satisfies that $i \in \ell(V(\mathcal{G}))$ whenever $0 \leq i<j$ and $j \in \ell(V(\mathcal{G}))$.
Thus the set $V(\mathcal{G})$ is finite or countable.
Moreover, it also implies that two vertices $s_1,s_2$ are adjacent if $\ell(s_1)=\ell(s_2) \pm 1$, while by definition of $\ell$, these are not adjacent if $\ell(s_1) \neq \ell(s_2) \pm 1$.
Thus the claim holds.\\
{\bf (ii)} First, take a vertex $t$ of $\mathcal{G}'$ and a simple path $P$ in $\mathcal{G}$ from $s_0$ to $t$.
Then the vertex $s$ of $P$ next to $s_0$ is adjacent to $s_0$ and contained in $\mathcal{G}'$.
On the other hand, if $\mathcal{G}'$ contains two vertices adjacent to $s_0$, then $s_0$ and a path in $\mathcal{G}'$ between these two vertices form a closed path in $\mathcal{G}$.
This is a contradiction, so that the first claim follows.
Since $s$ has degree $\leq 2$ in $\mathcal{G}$ and adjacent to $s_0 \not\in V(\mathcal{G}')$, $s$ is a terminal vertex of $\mathcal{G}'$.
Now the second claim is deduced by applying (i) to $\mathcal{G}'$ and $s$.\\
{\bf (iii)} This follows from (ii), since $\mathcal{G}$ is nonempty and has no terminal vertices.
\end{proof}
\begin{prop}
\label{prop:nofinitecomponent}
Let $(W,S)$ be an irreducible Coxeter system of an arbitrary rank, with Coxeter graph $\varGamma$.
Suppose that $|W|=\infty$ and $|W_I|<\infty$ for all finite subsets $I \subset S$.
Then $\varGamma \simeq \varGamma(A_{\infty})$, $\varGamma(B_{\infty})$, $\varGamma(D_{\infty})$ or $\varGamma(A_{\infty,\infty})$.
\end{prop}
\begin{proof}
In this proof, a full subgraph $\varGamma_I$ of $\varGamma$ is said to be {\it forbidden} if $|I|<\infty$ and $|W_I|=\infty$.
The hypothesis means that $|W|=\infty$ and $\varGamma$ is connected and contains no forbidden subgraphs.
This implies $|S|=\infty$ immediately.\\ \indent
{\bf Step 1. $\varGamma$ is acyclic:}\ 
This follows immediately from the fact that any nontrivial cycle in $\varGamma$ forms a forbidden subgraph.\\ \indent
{\bf Step 2. No $s \in S$ has degree $\geq 4$ in $\varGamma$:}\ 
Otherwise, this $s$ and the four adjacent vertices form a forbidden subgraph of $\varGamma$.
This is a contradiction.\\ \indent
{\bf Step 3. At most one $s \in S$ has degree $3$ in $\varGamma$:}\ 
Assume contrary that two distinct vertices $s,t \in S$ have degree $3$.
Since $\varGamma$ is connected, there is a path $P$ in $\varGamma$ between $s$ and $t$.
Then $s$, $t$, $P$ and all the vertices adjacent to $s$ or $t$ form a forbidden subgraph.
This is a contradiction.\\ \indent
{\bf Step 4. If some $s \in S$ has degree $3$ in $\varGamma$, then $\varGamma \simeq \varGamma(D_\infty)$:}\ 
By Steps 1--3, we can apply Lemma \ref{lem:lemmafortree} (ii) to this case.
This lemma shows that $\varGamma_{S \smallsetminus \{s\}}$ consists of three connected components $\simeq \varGamma(A_{n_1})$, $\varGamma(A_{n_2})$, $\varGamma(A_{n_3})$ (as unlabelled graphs) respectively, of which a terminal vertex is adjacent to $s$ in $\varGamma$.
By symmetry, we may assume $n_1 \geq n_2 \geq n_3 \geq 1$.\\ \indent
Now we have $n_1=\infty$ since $|S|=\infty$.
If $n_2 \geq 2$, then $\varGamma$ must contain a forbidden subgraph ($\simeq \varGamma(\widetilde{E_8})$ as unlabelled graphs), but this is a contradiction.
Thus we have $n_2=n_3=1$ and so $\varGamma \simeq \varGamma(D_\infty)$ as unlabelled graphs.
Moreover, every edge of $\varGamma$ must have no label (or label `$3$'), since otherwise $\varGamma$ must contain a forbidden subgraph again.
Hence $\varGamma \simeq \varGamma(D_\infty)$ (as Coxeter graphs) in this case.\\ \indent
{\bf Step 5. If all vertices of $\varGamma$ have degree $\leq 2$, then $\varGamma \simeq \varGamma(A_\infty)$, $\varGamma(B_\infty)$ or $\varGamma(A_{\infty,\infty})$:}\ 
First, we consider the case that $\varGamma$ has a terminal vertex.
Then Lemma \ref{lem:lemmafortree} (i) implies that $\varGamma \simeq \varGamma(A_\infty)$ as unlabelled graphs (note that $|S|=\infty$).
Moreover, by a similar argument to Step 4, the hypothesis ($\varGamma$ contains no forbidden subgraphs) detects the edge-labels of $\varGamma$, so that we have $\varGamma \simeq \varGamma(A_\infty)$ or $\varGamma(B_\infty)$ (as Coxeter graphs).
The other case is similar; we have $\varGamma \simeq \varGamma(A_{\infty,\infty})$ as Coxeter graphs by Lemma \ref{lem:lemmafortree} (iii) and the hypothesis.
Hence the proof is concluded.
\end{proof}
\begin{proof}[Proof of Theorem \ref{thm:coreofnormalizer} (for infinite $W$)]
Note that $Z(W)=1$ in this case.\\ \indent
{\bf Case 1. $(W,S)=(W(\mathcal{T}_n),S(\mathcal{T}_n))$ for $\mathcal{T}_n=A_\infty$, $B_\infty$, $D_\infty$ or $A_{\infty,\infty}$:}\ 
Put $L=1$ if $\mathcal{T}_n=B_\infty$, $L=2$ if $\mathcal{T}_n=D_\infty$.
Moreover, for $k \geq 1$, put
\begin{equation*}
J_k=\{s_1,s_2,\dotsc ,s_k\} \text{ if } \mathcal{T}_n \neq A_{\infty,\infty},\ J_k=\{s_{-k},s_{-k+1},\dotsc ,s_k\} \text{ if } \mathcal{T}_n=A_{\infty,\infty}.
\end{equation*}
\indent
{\bf Subcase 1-1. $\mathcal{T}_n=B_\infty$ or $D_\infty$, and $I=S(\mathcal{T}_k)$ for some $L \leq k<\infty$:}\ 
This is a case (i) or (ii) (for $\tau$ identity), so that we have to show $C_I=G_{\mathcal{T}_\infty}$.
Put $G_i=W_{J_{k+i}}$ and $H_i=N_{G_i}(W_I)$ for $i \geq 1$.
Then we have $\bigcup_{i=1}^{\infty}G_i=W$ and $\bigcup_{i=1}^{\infty}H_i=N_W(W_I)$, so that $C_I \subset \bigcup_{i=1}^{\infty}{\rm Core}_{G_i}(H_i)$ by Lemma \ref{lem:coreoflimit}.
Moreover, by the result of finite case (Section \ref{sec:coreoffinitecase}), we have ${\rm Core}_{G_i}(H_i)=G_{\mathcal{T}_{k+i}}$ for all $i \geq 1$.
Since $\bigcup_{i=1}^{\infty}G_{\mathcal{T}_{k+i}}=G_{\mathcal{T}_\infty}$ (cf.\ Lemma \ref{lem:defofpropersubgroupG}), we have $C_I \subset G_{\mathcal{T}_\infty}$.\\ \indent
On the other hand, we have $C_{S(\mathcal{T}_L)} \overset{\text{EL}}{\subset} C_I$, while $G_{\mathcal{T}_\infty} \subset Z_W(w_0(S(\mathcal{T}_L)))$ since $w_0(S(\mathcal{T}_L)) \in G_{\mathcal{T}_\infty}$ and $G_{\mathcal{T}_\infty}$ is abelian.
Thus $G_{\mathcal{T}_\infty} \subset N_W(W_{S(\mathcal{T}_L)})$ by Lemma \ref{lem:Zoflongestelement}, $G_{\mathcal{T}_\infty} \subset C_{S(\mathcal{T}_L)}$ by (\ref{eq:coreinsubgroup}) and so $G_{\mathcal{T}_\infty} \subset C_I$.
Hence $C_I=G_{\mathcal{T}_\infty}$.\\ \indent
{\bf Subcase 1-2. The hypothesis of Subcase 1-1 is not satisfied:}\ 
This is not a case (i) or (ii), so that we have to show $C_I=1$.\\ \indent
{\bf Subsubcase 1-2-1. $|I|<\infty$:}\ 
Let $w \in C_I$.
Now take a sufficiently large $4 \leq k<\infty$ so that $I \subset J_k$ and $w \in W_{J_k}$.
Put $G_i=W_{J_{k+i}}$ and $H_i=N_{G_i}(W_I)$ for $i \geq 1$, so that $\bigcup_{i=1}^{\infty}G_i=W$ and $\bigcup_{i=1}^{\infty}H_i=N_W(W_I)$.
Now by the hypothesis of Subcase 1-2, and by the result for finite case (Section \ref{sec:coreoffinitecase}), we have ${\rm Core}_{G_i}(H_i) \subset Z(G_i) \subset \{1,w_0(J_{k+i})\}$ for all $i$.
Moreover, by Lemma \ref{lem:coreoflimit}, we have $C_I \subset \bigcup_{i=1}^{\infty}{\rm Core}_{G_i}(H_i)$.
Since $w_0(J_{k+i}) \not\in W_{J_k}$ for any $i \geq 1$, this implies that $w=1$ by the choice of $k$.
Hence we have $C_I=1$.\\ \indent
{\bf Subsubcase 1-2-2. $|I|=\infty$:}\ 
If $I$ has an irreducible component $J$ of infinite cardinality, then $C_I=1$ by (\ref{eq:coreofinfinite}).
Thus we may assume that $I$ is a union of infinitely many irreducible components of finite cardinality.
Now we can choose indices $4 \leq i \leq j<\infty$ so that $s_k \not\in I$ for all $i \leq k \leq j$, $s_{i-1} \in I$ and $s_{j+1} \in I$.
Let $K_1$, $K_2$ be the (distinct) irreducible components of $I$ containing $s_{i-1}$, $s_{j+1}$ respectively.
Then we have $C_I \overset{\text{CL}}{\subset} C_{I \smallsetminus (K_1 \cup K_2)}$ and so $C_I \subset C_{K_1 \cup K_2}$ by (\ref{eq:coreofseparatedI}).
Moreover, we have $C_{K_1 \cup K_2}=1$ by Subsubcase 1-2-1.
Thus $C_I=1$.\\ \indent
{\bf Case 2. $(W,S) \not\simeq (W(\mathcal{T}),S(\mathcal{T}))$ for $\mathcal{T}=A_\infty$, $B_\infty$, $D_\infty$, $A_{\infty,\infty}$:}\ 
This is not a case (i) or (ii), so that we have to show $C_I=1$.
By Proposition \ref{prop:nofinitecomponent}, there is a finite subset $J_0 \subset S$ such that $|W_{J_0}|=\infty$.
This $J_0$ consists of only finitely many irreducible components, and so we have $|W_J|=\infty$ for some irreducible component of $J_0$.
Since $\varGamma$ is connected and $|J|<\infty$, there is a (finite) sequence $s_1,s_2,\dotsc ,s_r$ of elements of $S$ such that $s_i \not\in I_{i-1} \cup {I_{i-1}}^\perp$ for all $1 \leq i \leq r$ and $J \subset I_r$, where we put $I_0=I$ and $I_i=I_{i-1} \cup \{s_i\}$ ($1 \leq i \leq r$) inductively.
Now we have $C_{I_{i-1}} \overset{\text{EL}}{\subset} C_{I_i}$ for all $1 \leq i \leq r$, so that $C_I \subset C_{I_{r-1}}$ and $C_I \subset C_{I_r}$.\\ \indent
{\bf Subcase 2-1. $I_r \neq S$:}\ 
Now an irreducible component of $I_r$ (namely, the one containing $J$) generates an infinite group.
Thus $C_I \subset C_{I_r}=1$ by (\ref{eq:coreofinfinite}).\\ \indent
{\bf Subcase 2-2. $I_r=S$:}\ 
Note that $r \geq 1$ since $I$ is proper.
Since $(W,S)$ is irreducible, we have $C_I \subset C_{I_{r-1}}=1$ by (\ref{eq:coreofmaximalparabolic}).
Hence the proof is concluded.
\end{proof}
%

\ \\
{\bf Koji NUIDA}\medskip\\
Graduate School of Mathematical Sciences, University of Tokyo\\
3-8-1 Komaba, Meguro-ku, Tokyo, 153-8914, Japan\medskip\\
E-mail: nuida@ms.u-tokyo.ac.jp
\end{document}